\documentclass[11pt, a4paper]{article}
\usepackage{amsmath,amssymb,latexsym,color}
\usepackage{authblk}
\usepackage[mathscr]{eucal}
\usepackage[colorlinks,
            linkcolor=blue,
            anchorcolor=green,
            citecolor=magenta
           ]{hyperref}
\newtheorem{thm}{Theorem}[section]

\newtheorem{lemma}[thm]{Lemma}

\newtheorem{hyp}{Hypothesis}
\newtheorem{rem}[thm]{Remark}

\newtheorem{example}{Example}[section]
\newtheorem{defin}{Definition}[section]

\newcommand{\proof}{{\it Proof.\quad}}
\newcommand{\qed}{\hfill\Box\medskip}

\usepackage{CJK}

\begin{document}
%\begin{CJK*}{GBK}{song}

\renewcommand{\baselinestretch}{1.3}
%%%%%%%%%%%%%%%%%%%%%%%%%%%%%%%%%%%%%%%%%%%%%%%%%%%%%%%%%%%%%%%%%%%%%%%%%%%%%%%%%%%%%%%%
%%%%%%%%%%%%%%%%%%%%%%%%%%%%%%%%%%%%%%%%%%%%%%%%%%%%%%%%%%%%%%%%%%%%%%%%%%%%%%%%%%%%%%%%
\title{\bf Non-trivial $t$-intersecting families for the distance-regular graphs of bilinear forms}

\author[1,2]{Mengyu Cao\thanks{E-mail: \texttt{caomengyu@mail.bnu.edu.cn}}}
\author[2]{Benjian Lv\thanks{Corresponding author. E-mail: \texttt{bjlv@bnu.edu.cn}}}
\author[2]{Kaishun Wang\thanks{E-mail: \texttt{wangks@bnu.edu.cn}}}
\affil[1]{\small Department of Mathematical Sciences, Tsinghua University, Beijing 100084, China}
\affil[2]{\small Laboratory of Mathematics and Complex Systems (Ministry of Education), School of Mathematical Sciences, Beijing Normal University, Beijing 100875, China}
 \date{}
 \maketitle

\begin{abstract}

Let $V$ be an $(n+\ell)$-dimensional vector space over a finite field, and $W$ a fixed $\ell$-dimensional subspace of $V$. Write ${V\brack n,0}$ to be the set of all $n$-dimensional subspaces $U$ of $V$ satisfying $\dim(U\cap W)=0$.  A family $\mathcal{F}\subseteq{V\brack n,0}$ is $t$-intersecting if $\dim(A\cap B)\geq t$ for all $A,B\in\mathcal{F}$. A $t$-intersecting family $\mathcal{F}\subseteq{V\brack n,0}$ is called non-trivial if $\dim(\cap_{F\in\mathcal{F}}F)<t$. In this paper, we describe the structure of non-trivial $t$-intersecting families of ${V\brack n,0}$ with large size. In particular, we show the structure of the non-trivial $t$-intersecting families with maximum size, which extends the Hilton-Milner Theorem for ${V\brack n,0}$.



\end{abstract}
\section{Introduction}

Let $X$ be a set of size $n$, and ${X\choose k}$ denote the set of all $k$-subsets of $X$. For a positive integer $t$, the family $\mathcal{F}\subseteq {X\choose k}$ is said to be \emph{$t$-intersecting} if $|A \cap B|\geq t$ for all $A, B\in\mathcal{F}$. A $t$-intersecting family is called \emph{trivial} if all its members contain a common specified $t$-subset of $X$, and \emph{non-trivial} otherwise.

The famous Erd\H{o}s-Ko-Rado Theorem \cite{Erdos-Ko-Rado-1961-313,Frankl-1978,Wilson-1984} showed that each $t$-intersecting family of ${X\choose k}$ with maximum size is a trivial family consisting of all $k$-subsets that contain a fixed $t$-subset of $X$ for $n>(t +1)(k-t +1)$. In \cite{Frankl-1978}, Frankl  made a conjecture on the maximum size of a $t$-intersecting family of ${X\choose k}$ for all positive integers $t,\ k$ and $n$. This conjecture was partially proved by Frankl and F\"{u}redi  \cite{Frankl--Furedi-1991} and completely settled by Ahlswede and Khachatrian \cite{Ahlswede-Khachatrian-1997}. Determining the structure of non-trivial $t$-intersecting families of ${X\choose k}$ with maximum size was a long-standing problem. The first result was the Hilton-Milner Theorem \cite{Hilton-Milner-1967, Frankl-Furedi-1986} which describes the structure of such families for $t=1$. A significant step was taken in \cite{Frankl-1978-1} by Frankl, who determined such families for $t\geq 2$ and $n>n_1(k,t)$.  Ahlswede and Khachatrian \cite{Ahlswede-Khachatrian-1996} gave a complete result on non-trivial intersection problems for finite sets.

Recently, other non-trivial $t$-intersecting families had been studied. Han and Kohayakawa \cite{Han-Kohayakawa} determined the structure of the third largest maximal $1$-intersecting families of ${X\choose k}$ with $6\leq 2k<n$.  Kostochka and Mubayi \cite{Kostochka-Mubayi} described the structure of $1$-intersecting families of ${X\choose k}$ for large $n$ whose size is quite a bit smaller than ${n-1\choose k-1}$. In \cite{Cao-set}, we described the structure of maximal non-trivial $t$-intersecting families of ${X\choose k}$ with large size.

The Erd\H{o}s-Ko-Rado Theorem and the Hilton-Milner Theorem for finite sets have natural extensions to vector spaces. Let $q$ be a prime power. Recall that  for  positive integers $a$ and $b$ the \emph{Gaussian binomial coefficient} is defined by
\begin{align*}
{a\brack b} = \prod_{0\leq i<b}\frac{q^{a-i}-1}{q^{b-i}-1}.
\end{align*}
In addition, we set ${a\brack 0}=1$ and ${a\brack c} =0$ if $c$ is a negative integer, and denote ${a\brack 1}$ by $\theta_a$.

Let $V$ be a $v$-dimensional vector space over the finite field $\mathbb{F}_q$,  and ${V\brack k}$ denote the set of all $k$-dimensional subspaces (`$k$-subspaces' for short) of $V$. Note that the size of ${V\brack k}$ is  ${v\brack k}$.  For positive integer $t$, a family $\mathcal{F}\subseteq{V\brack k}$ is called $t$-\emph{intersecting} if $\dim(A\cap B)\geq t$ for all $A,B\in\mathcal{F}$. A $t$-intersecting family $\mathcal{F}\subseteq{V\brack k}$ is called \emph{trivial} if all its members contain a common $t$-subspace of $V$, and \emph{non-trivial} otherwise. The structure of $t$-intersecting families of ${V\brack k}$ with maximum size were well studied \cite{Deza-Frankl-1983,PR,Hsieh-1975-1,Tanaka-2006-903}, which are known as the Erd\H{o}s-Ko-Rado Theorem for vector spaces. Using the covering number, Blokhuis et al. \cite{AB} obtained a vector space version of the Hilton-Milner Theorem, which describes the structure of non-trivial $1$-intersecting families of ${V\brack k}$ with maximum size. In \cite{Cao-vec}, we described the structure of maximal non-trivial $t$-intersecting families of ${V\brack k}$ with large size, from which the non-trivial $t$-intersecting families with maximum size are determined (Also see \cite{D'haeseleer}).

Let $V$ be an $(n+\ell)$-dimensional row vector space over $\mathbb{F}_{q}$ in the following, and $W$ a fixed $\ell$-subspace of $V$. For integers $m$ and $s$, we say that a subspace $P$ has {\em type} $(m, s)$ or $P$ is an $(m,s)$-subspace of $V$ if $\dim (P)
= m$ and $\dim(P\cap W)=s.$ For each subspace $S$ of $V$, let ${S\brack m,s}$ denote the set of all $(m,s)$-subspaces of $V$ which are contained in $S$. Note that the set ${V\brack m,s}$ forms an orbit under the action of the singular general linear group $GL_{n+\ell,n}(\mathbb{F}_q)$ on the set of all subspaces of $V$ when ${V\brack m,s}\neq\emptyset$ (\cite{KWW}), and the configuration $({V\brack m,s},\ \Lambda)$ forms an association scheme, where $\Lambda$ is the set of all orbits of the action of $GL_{n+\ell,n}(\mathbb{F}_q)$ on ${V\brack m,s}\times{V\brack m,s}$. It is well known that $({V\brack n,0},\ \Lambda)$ is the \emph{bilinear forms scheme}, and the \emph{bilinear forms graph} $H_{q}(\ell,n)$ is a distance regular graph defined on ${V\brack n,0}$, with two vertices $A$ and $B$ adjacent if $\dim(A\cap B)=n-1.$ A bilinear forms graph usually can be viewed as a subgraph of a Grassmann graph, or a $q$-analogue of a Hamming graph. Recall that $H_q(\ell,n)$ and $H_q(n,\ell)$ are isomorphic, and the distance between two vertices $A$ and $B$ is $i$ in $H_q(\ell,n)$ with $\ell\geq n$ if and only if $\dim(A\cap B)=n-i$. See \cite{ABE} for more information about $H_{q}(\ell,n).$

A family $\mathcal{F}\subseteq{V\brack m,k}$ is
\emph{$t$-intersecting} if $\dim(A\cap B)\geq t$  for all $A,
B\in\mathcal{F}.$ A $t$-intersecting family $\mathcal{F}\subseteq{V\brack m,k}$ is called \emph{trivial} if all the elements in  $\mathcal{F}$ contains a common $t$-subspace  and \emph{non-trivial} otherwise, and is called \emph{maximal} if $\mathcal{F}\cup\{A\}$ is not $t$-intersecting for each $A\in{V\brack n,0}\setminus\mathcal{F}$. It is obvious that for $\ell\geq n$, $\mathcal{F}\subseteq {V\brack n,0}$ is $t$-intersecting if and only if $\mathcal{F}$ is a subset of vertices with maximum distance $n-t$ in $H_{q}(\ell,n).$ In \cite{Huang,Tanaka-2006-903}, the authors proved the Erd\H{o}s-Ko-Rado Theorem
for ${V\brack n,0}$ using different methods, and Ou et al. \cite{LO} generalized this result to  ${V\brack m,k}.$ In \cite{Gong} and \cite{Hou}, the authors proved the Hilton-Milner Theorem for ${V\brack n,0}$ and ${V\brack m,0}$, respectively.

In this paper, we consider the structure of maximal non-trivial $t$-intersecting families of ${V\brack n,0}$ with $\ell\geq n$. When $t=n-1$, from the structure of maximal $(n-1)$-intersecting families of ${V\brack n}$ or \cite[Lemma~17]{Hemmeter}, one can deduces that each maximal $(n-1)$-intersecting family of ${V\brack n,0}$ is a collection of all $(n,0)$-subspaces of $V$ containing some fixed $(n-1,0)$-subspace of $V$, or a collection of all $(n,0)$-subspaces of $V$ contained in some $(n+1,1)$-subspace of $V$. Thus, we only need to consider the case when $3\leq n\leq \ell$ and $1\leq t\leq n-2$.

The key notion used in this paper is the `$t$-covering number' of a $t$-intersecting family, which is a generalization of the classical `covering number'. We refer the readers to \cite{Cao-set, Cao-vec, Frankl-1980, Frankl--Ota--Tokushige-1991, Furedi-1988,Furuya--Takatou-1988} for more results about the covering number and $t$-covering number. The $t$-{\em covering number} $\tau_t(\mathcal{F})$ of a family $\mathcal{F}\subseteq {V\brack m,0},$ is the minimum dimension of a subspace $T$ of $V$ such that $T\cap W=0$ and $\dim(T\cap F)\geq t$ for all $F\in \mathcal{F}$.  Let $\mathcal{F}\subseteq{V\brack m,0}$ be a $t$-intersecting family. Observe that $t\leq \tau_t(\mathcal{F})\leq m,$ and $\mathcal{F}$ is trivial if and only if $\tau_t(\mathcal{F})=t.$ To present our results let us first introduce the following constructions of non-trivial $t$-intersecting families.

\

\noindent\textbf{Family I.}\quad Let $X\in{V\brack t,0}$ and $M\in{V\brack n+1,1}$ be two subspaces with $X\subseteq M$. Denote
\begin{align*}
\mathcal{H}_1(X,M)=\left\{F\in{V\brack n,0} \mid X\subseteq F,\ \dim(F\cap M)\geq t+1\right\}\cup {M\brack n,0}.
\end{align*}

\noindent\textbf{Family II.}\quad Let $X\in{V\brack t,0},$ $M\in{V\brack n,k}$ and $C\in{V\brack c,c-n}$ be three subspaces with $X\subseteq M\subseteq C,$ where $k\in\{0,1\}$ and $c\in\{n+1, n+2,\cdots, 2n-t, n+\ell\}$.  Denote
\begin{align*}
	\mathcal{H}_2(X,M,C)=\mathcal{A}(X,M,C)\cup\mathcal{B}(X,M,C)\cup\mathcal{C}(X,M,C),\quad \mbox{where}
\end{align*}
\begin{align*}
\mathcal{A}(X,M,C)=&\left\{F\in{V\brack n,0}\mid X\subseteq F,\ \dim(F\cap M)\geq t+1\right\},\\
\mathcal{B}(X,M,C)=&\left\{F\in{V\brack n,0}\mid F\cap M=X,\ \dim(F\cap C)=c-n+t\right\},\\
\mathcal{C}(X,M,C)=&\left\{F\in{C\brack n,0}\mid \dim(F\cap X)=t-1,\ \dim(F\cap M)=n-1\right\}.
\end{align*}

\noindent\textbf{Family III.}\quad Let $Z$ be a $(t+2,k)$-subspace of $V$, where $k\in\{0,1\}$. Denote
\begin{align*}
\mathcal{H}_3(Z)=\left\{F\in{V\brack n,0}\mid \dim(F\cap Z)\geq t+1\right\}.
\end{align*}

\begin{rem}\label{b-obs}
{\rm In Family~II, if $\dim(C)=n+1$, then $\mathcal{H}_2(X,M,C)=\mathcal{H}_1(X,C)$; if $C=V$, then $\mathcal{H}_2(X,M,V)=\mathcal{A}(X,M,V)\cup\mathcal{C}(X,M,V)$; if $t=n-2$, then $\mathcal{H}_2(X,M,V)=\mathcal{H}_3(M)$.}
\end{rem}

 Our first main result describes the structure of all maximal non-trivial $t$-intersecting families of $(n,0)$-subspaces of $V$ with large size.
\begin{thm}\label{b-main1}
Let $1\leq t\leq n-2$, $n+t+2\leq \ell$, and $(\ell,q)\neq(n+t+2,2)$ or $(n+t+3,2)$. If $\mathcal{F}\subseteq{V\brack n,0}$ is a maximal non-trivial $t$-intersecting family and
\begin{align*}
|\mathcal{F}|\geq q^{\ell(n-t-1)+1}\theta_{n-t-1}-q^{\ell(n-t-2)+3}{n-t-1\brack 2},
\end{align*}
then one of the following holds:
\begin{itemize}
\item[{\rm (i)}]
$\mathcal{F}=\mathcal{H}_2(X,M,C)$ for some $(t,0)$-subspace $X$, $(n,k)$-subspace $M$ and $(c,c-n)$-subspace $C$ of $V$ with $X\subseteq M\subseteq C$, where $k\in\{0,1\}$ and $c\in\{n+1,n+2,\ldots, 2n-t, n+\ell\}$;
\item[{\rm (ii)}] $\mathcal{F}=\mathcal{H}_3(Z)$ for some $(t+2,0)$-subspace $Z$ of $V$, and $\frac{n}{2}-1\leq t\leq n-2$;
\item[{\rm (iii)}]  $\mathcal{F}=\mathcal{H}_3(Z)$ for some $(t+2,1)$-subspace $Z$ of $V$, and $\frac{n}{2}-\frac{1}{2}\leq t\leq n-2$.
\end{itemize}
\end{thm}

By comparing the size of each family in Theorem~\ref{b-main1}, we can obtain the structure of non-trivial $t$-intersecting families with maximum size, which extends the Hilton-Milner Theorem for $(n,0)$-subspaces of $V$ given in \cite{Gong}.
\begin{thm}\label{b-main2}
Let  $1\leq t\leq n-2$, $n+t+2\leq \ell$, and $(\ell,q)\neq(n+t+2,2)$ or $(n+t+3,2)$. Let $\mathcal{F}\subseteq{V\brack n,0}$ be a non-trivial $t$-intersecting family. Then the following hold.
\begin{itemize}
\item[{\rm(i)}] If $1\leq t\leq \frac{n}{2}-1,$ then $|\mathcal{F}|\leq q^{\ell(n-t)}-\prod_{j=1}^{n-t}(q^\ell-q^j)+q^{n-t}(q^t-1)$. Equality holds if and only if $\mathcal{F}=\mathcal{H}_1(X,M)$ for some $(t,0)$-subspace $X$ and $(n+1,1)$-subspace $M$ of $V$ with $X\subseteq M$.
\item[{\rm(ii)}] If $\frac{n}{2}-\frac{1}{2}\leq t\leq n-2$, then $|\mathcal{F}|\leq q^{\ell(n-t-1)}\theta_{t+2}-q^{\ell(n-t-2)+1}\theta_{t+1}$. Equality holds if and only if $\mathcal{F}=\mathcal{H}_3(Z)$ for some $(t+2,0)$-subspace $Z$ of $V$.
\end{itemize}
\end{thm}

In Section 2, we will give some upper bounds for the non-trivial $t$-intersecting families of ${V\brack n,0}$. In Sections 3 and 4, we will prove Theorems~\ref{b-main1} and~\ref{b-main2}, respectively.
\section{Upper bounds for non-trivial $t$-intersecting families}
In this section, we give some upper bounds for the maximal non-trivial $t$-intersecting families of $(n,0)$-subspaces of $V$. Firstly, we state several useful results about the number of subspaces in a given family. It is well known that for positive integers  $m$ and $i$ with $i<m$,
\begin{align}
	{m\brack i}={m-1\brack i-1}+q^i{m-1\brack i},\label{b0}\\
	q^{m-i}<\frac{q^m-1}{q^i-1}<q^{m-i+1}.\label{b00}
\end{align}
\iffalse
\begin{lemma}{\rm(\cite[Lemma~9.3.2]{ABE})}\label{lem4}
	Suppose $0\leq i,j\leq n.$ If $X$ is a $j$-subspace of $V$, then there are precisely $q^{(i-m)(j-m)}{n-j\brack i-m}{j\brack m}$ $i$-subspaces $Y$ in $V$ such that $X\cap Y$ is an $m$-subspace.
\end{lemma}
\fi

Let $N'(m_{1},h_{1};m,h;n+\ell,n)$ be the number of $(m,h)$-subspaces in $V$ containing a given $(m_{1},h_{1})$-subspace of $V$. The next lemma states the exact value of $N'(m_{1},h_{1};m,h;n+\ell,n)$.
\begin{lemma}{\rm{(\cite{KW})}}\label{lem5}
Let $0\leq h_{1}\leq h\leq \ell$ and $0\leq m_{1}-h_{1}\leq m-h\leq n$. Then,
	$${\textstyle N^{'}(m_{1},h_{1};m,h;n+\ell,n)=q^{(\ell-h)(m-h-m_{1}+h_{1})}{{n-(m_{1}-h_{1})}\brack{(m-h)-(m_{1}-h_{1})}}{{\ell-h_{1}}\brack{h-h_{1}}}.}$$
\end{lemma}

For $\mathcal{F}\subseteq {V\brack n,0}$ and a subspace $S$ of $V$ with $S\cap W=0$, let $\mathcal{F}_S$ denote the collection of all $(n,0)$-subspaces in $\mathcal{F}$ which contain $S$. For convenience, we set $\prod_{i=a}^{a-1}b_i=1$ in the following.

\begin{lemma}\label{b-lem7}
	Let $\mathcal{F}\subseteq{V\brack n,0}$ be a $t$-intersecting family and $S$ an $(s,0)$-subspace of $V$ with $t-1\leq s\leq n-1.$ If there is an  $F^\prime\in \mathcal{F}$ such that $\dim(S\cap F^\prime)=r\leq t-1,$ then $n\geq s+t-r$, and for each $i\in\{1,2,\ldots,t-r\}$ there exists an $(s+i,0)$-subspace $T_i$ with $S\subseteq T_i$ such that $|\mathcal{F}_S|\leq q^{i(s-r)}{n-s\brack i}|\mathcal{F}_{T_i}|$.
\end{lemma}
\proof From $\dim(S\cap F^{\prime})=r,$ observe that the type of $S+F^\prime$ is  $(n+s-r,s-r)$. For every $i\in\{1,2,\ldots,t-r\}$, set
\begin{align*}
	\mathcal{H}_i=\left\{H\in{S+F^\prime\brack s+i,0}\mid S\subseteq H\right\}.
\end{align*}
For each $F\in\mathcal{F}_S,$ since $\dim(F\cap F^\prime)\geq t$, we have $\dim(F\cap (S+F^\prime))\geq s+t-r$. Then $n\geq s+t-r$, and
 there exists $H\in\mathcal{H}_i$ such that $H\subseteq F,$ implying that $\mathcal{F}_S=\cup_{H\in\mathcal{H}_i}\mathcal{F}_H.$ Let $T_i$ be an element in $\mathcal{H}_i$ satisfying $|\mathcal{F}_H|\leq|\mathcal{F}_{T_i}|$ for each $H\in\mathcal{H}_i.$ Since $|\mathcal{H}_i|=N^\prime(s,0;s+i,0;n+s-r,n)=q^{i(s-r)}{n-s\brack i}$ by Lemma~\ref{lem5}, we obtain the required result.  $\qed$

Note that $|\mathcal{F}_{T_i}|\leq N^\prime(s+i,0;n,0;n+\ell,n)=q^{\ell(n-s-i)}$ for each $(s+i,0)$-subspace $T_i$ of $V$ by Lemma~\ref{lem5}. Then we can obtain the following lemma by setting $i=t-r$ in Lemma~\ref{b-lem7}.
\begin{lemma}\label{b-lem8}
	Let $\mathcal{F}\subseteq{V\brack n,0}$ be a $t$-intersecting family and $S$ an $(s,0)$-subspace of $V$ with $t-1\leq s\leq n-1.$ If there exists $F^\prime\in \mathcal{F}$ such that $\dim(S\cap F^\prime)=r\leq t-1,$ then $|\mathcal{F}_S|\leq q^{\ell(n-s-t+r)+(t-r)(s-r)}{n-s\brack t-r}$.
\end{lemma}

\begin{lemma}\label{b-prop7}
	Let $1\leq t\leq n-2\leq\ell-2$, and $\mathcal{F}\subseteq {V\brack n,0}$ a maximal $t$-intersecting family with $t+2\leq \tau_t(\mathcal{F})=m\leq n$. Then
\begin{align}\label{ll-2}
	|\mathcal{F}|\leq q^{\ell(n-m)+2m-2t-1}\theta_{n-m+1}\theta_{n-m+2}{m\brack t}\prod_{j=1}^{m-t-2}(q^{t+j-1}\theta_{n-t-j+1}).
\end{align}
	Moreover, if $\ell\geq n+t+1,$ then
\begin{align}\label{ll-3}
	|\mathcal{F}|\leq q^{\ell(n-t-2)+3}\theta_{n-t-1}\theta_{n-t}{t+2\brack 2}.
\end{align}
\end{lemma}
\proof Let $T$ be an $(m,0)$-subspace of $V$ which satisfies $\dim(T\cap F)\geq t$ for all $F\in\mathcal{F}$. Then $\mathcal{F}=\bigcup_{H\in{T\brack t,0}}\mathcal{F}_H,$ and hence $
|\mathcal{F}|\leq{m\brack t}|\mathcal{F}_{H_1}|$ for some $H_1\in{T\brack t}$.
We claim that there exists $H^\prime\in {V\brack m-2,0}$ such that
\begin{align}\label{ll-1}
	|\mathcal{F}|\leq {m\brack t}\left(\prod_{j=1}^{m-t-2}q^{t+j-1}\theta_{n-t-j+1}\right)|\mathcal{F}_{H^\prime}|
\end{align}
Actually, if $m=t+2$, then it is clear that (\ref{ll-1}) holds; if $m\geq t+3$, using Lemma~\ref{b-lem7} repeatedly, there exist $H_2\in{V\brack t+1,0}$, $H_3\in{V\brack t+2,0},\ldots,H_{m-t-1}\in{V\brack m-2,0}$ such that $H_j\subseteq H_{j+1}$ and $|\mathcal{F}_{H_j}|\leq q^{t+j-1}\theta_{n-t-j+1}|\mathcal{F}_{H_{j+1}}|$ for every $j\in\{1,2,\ldots,m-t-2\}$, implying that (\ref{ll-1}) holds.

Since $\tau_t(\mathcal{F})>m-2$, there exists $F\in\mathcal{F}\setminus\mathcal{F}_{H^\prime}$ such that $\dim(F\cap H^\prime)\leq t-1$.
\
\medskip
\

\noindent\textbf{Case 1.}\quad There does not exist $F\in\mathcal{F}\setminus\mathcal{F}_{H^\prime}$ such that $\dim(F\cap H^\prime)=t-1$.

In this case, we have $t\geq 2$. Let $F_1$ be an element in $\mathcal{F}\setminus\mathcal{F}_{H^\prime}$ with $\dim(F_1\cap H^\prime)=a_1\leq t-2$. By Lemma~\ref{b-lem7}, note that $a_1\geq  \max\{0, m-2+t-n\}$. Let
\begin{align*}
g(a)=q^{\ell(n-m+2-t+a)+(t-a)(m-2-a)}{n-m+2\brack t-a}
\end{align*}
for $a\in\{\max\{0, m-2+t-n\},\ldots,t-2\}$. Then
\begin{align*}
\frac{g(a+1)}{g(a)}=\frac{q^{\ell-t+2a-m+3}(q^{t-a}-1)}{q^{n-m+3-t+a}-1}>q^{\ell-n+a}(q^{t-a}-1)>1
\end{align*}
for $a\in\{\max\{0, m-2+t-n\},\ldots,t-3\}$. That is, the function $g(a)$ is increasing as $s\in\{\max\{0, m-2+t-n\},\ldots,t-2\}$ increases. Hence, by Lemma~\ref{b-lem8} we have
\begin{align}\label{b-equ16}
|\mathcal{F}_{H^\prime}|\leq g(a_1)\leq g(t-2)=q^{\ell(n-m)+2(m-t)}{n-m+2\brack 2}.
\end{align}
\
\medskip
\

\noindent\textbf{Case 2.}\quad There exists $F_2\in\mathcal{F}\setminus\mathcal{F}_{H^\prime}$ such that $\dim(F_2\cap H^\prime)=t-1$.

By Lemma~\ref{b-lem7}, there exists an $(m-1,0)$-subspace $H^{\prime\prime}$ such that $H^\prime\subseteq H^{\prime\prime}$ and $|\mathcal{F}_{H^\prime}|\leq q^{m-1-t}\theta_{n-m+2}|\mathcal{F}_{H^{\prime\prime}}|.$
Since $\tau_t(\mathcal{F})> m-1$, there exists $F_3\in \mathcal{F}$ such that $\dim(F_3\cap H^{\prime\prime})\leq t-1.$

If $\dim(F_3\cap H^{\prime\prime})=t-1$, by Lemma~\ref{b-lem8}, then  $|\mathcal{F}_{H^{\prime\prime}}|\leq q^{\ell(n-m)+m-t}\theta_{n-m+1},$ implying that
\begin{eqnarray}\label{b-equ17}
	|\mathcal{F}_{H^\prime}|\leq q^{\ell(n-m)+2m-2t-1}\theta_{n-m+2}\theta_{n-m+1}.
\end{eqnarray}

Suppose $\dim(F_3\cap H^{\prime\prime})=a_2\leq t-2$. By Lemmas~\ref{b-lem7} and \ref{b-lem8}, note that $a_2\geq  \max\{0, m-1+t-n\}$ and
\begin{align*}
|\mathcal{F}_{H^{\prime\prime}}|\leq q^{\ell(n-m+1-t+a_2)+(t-a_2)(m-1-a_2)}{n-m+1\brack t-a_2}.
\end{align*}
Similar to Case 1, it is routine to verify that
\begin{align*}
|\mathcal{F}_{H^{\prime\prime}}|\leq q^{\ell(n-m-1)+2(m-t+1)}{n-m+1\brack 2}.
\end{align*}
Hence we have
\begin{eqnarray}\label{b-equ18}
	|\mathcal{F}_{H^\prime}|\leq q^{\ell(n-m-1)+3m-3t+1}\theta_{n-m+2}{n-m+1\brack 2}.
\end{eqnarray}

Since
\begin{align*}
\frac{q^{-1}\theta_{n-m+2}\theta_{n-m+1}}{{n-m+2\brack 2}}=q^{-1}(q+1)>1
\end{align*}
and
\begin{align*}
 \frac{q^{-1}\theta_{n-m+2}\theta_{n-m+1}}{q^{-\ell+m-t+1}\theta_{n-m+2}{n-m+1\brack 2}}=\frac{q^{\ell-m+t-2}(q^2-1)}{q^{n-m}-1}>q^{\ell-n+t-2}(q^2-1)\geq 1,
\end{align*}
we have (\ref{ll-2}) holds by (\ref{ll-1}--\ref{b-equ18}).

Let
\begin{align*}
p(m^\prime)=q^{\ell(n-m^\prime)+2m^\prime-2t-1}\theta_{n-m^\prime+1}\theta_{n-m^\prime+2}{m^\prime\brack t}\prod_{j=1}^{m^\prime-t-2}(q^{t+j-1}\theta_{n-t-j+1})
\end{align*}
for $m^\prime\in\{t+2,t+3,\ldots,n\}$. By $\ell\geq n+t+1$, we have
\begin{align*}
\frac{p(m^\prime)}{p(m^\prime+1)}=\frac{q^{\ell-m^\prime}(q^{m^\prime-t+1}-1)(q-1)}{(q^{m^\prime+1}-1)(q^{n-m^\prime}-1)}>\frac{q^{\ell-m^\prime}\cdot q^{m^\prime-t}}{q^{m^\prime+1}\cdot q^{n-m^\prime}}\geq 1
\end{align*}
for $m^\prime\in\{t+2,t+3,\ldots,n-1\}$. This yields that $|\mathcal{F}|\leq p(t+2)$, and hence (\ref{ll-3}) holds.    $\qed$

In the following, we consider maximal non-trivial $t$-intersecting families $\mathcal{F}\subseteq {V\brack n,0}$ with $\tau_t(\mathcal{F})=t+1.$
\begin{hyp}\label{b-hyp1}
	Suppose $4\leq n+1\leq \ell$ and $1\leq t\leq n-2$. Let $\mathcal{F}\subseteq {V\brack n,0}$ be a maximal non-trivial $t$-intersecting family with $\tau_t(\mathcal{F})=t+1.$ Define $\mathcal{T}$ to be the set of all $(t+1,0)$-subspaces $T$ of $V$ which satisfy $\dim(T\cap F)\geq t$ for all $F\in\mathcal{F}$. Set $M:=\sum_{T\in\mathcal{T}}T$.
\end{hyp}

\begin{lemma}\label{b-lem2}
	Let $n,\ \ell,\ t,\ \mathcal{F}$ and $\mathcal{T}$ be as in Hypothesis~\ref{b-hyp1}. If $3n-3t-2\leq \ell$, then $\mathcal{T}\subseteq {V\brack t+1,0}$ is a $t$-intersecting family with $t\leq\tau_t(\mathcal{T})\leq t+1$.
\end{lemma}
\proof It is clear that $t\leq\tau_t(\mathcal{T})\leq t+1$ holds if $\mathcal{T}$ is $t$-intersecting. For each $T\in \mathcal{T}$, by maximality, $\mathcal{F}$ contains all $(n,0)$-subspaces of $V$ containing $T$. In the following, we will prove that $\mathcal{T}$ is $t$-intersecting by reduction to absurdity. Suppose there exist $T_1,T_2\in\mathcal{T}$ such that $\dim(T_1\cap T_2)=r<t.$ It suffices to construct two $(n,0)$-subspaces $F_1$ and $F_2$ such that $T_1\subseteq F_1$, $T_2\subseteq F_2$ and $\dim(F_1\cap F_2)<t$. Assume that
\begin{align*}
T_1=\langle \beta_1,\ldots,\beta_r,\beta_{r+1},\ldots,\beta_{t+1}\rangle\quad\mbox{and}\quad T_2=\langle \beta_1,\ldots,\beta_r,\beta^\prime_{r+1},\ldots,\beta^\prime_{t+1}\rangle.
\end{align*}
Then $\beta_1,\ldots,\beta_r,\beta_{r+1},\ldots,\beta_{t+1},\beta^\prime_{r+1},\ldots,\beta^\prime_{t+1}$ are linearly independent.

By maximality of $\mathcal{F}$, there exist $A_1,A_2\in\mathcal{F}$ such that $T_1\subseteq A_1$ and $T_2\subseteq A_2$. Observe that $A_1+A_2$ is of type $(\dim(A_1+A_2),\dim(A_1+A_2)-n)$. Since $\dim(A_1\cap A_2)\geq t$, we have $\dim(A_1+A_2)\leq 2n-t\leq \ell-n+2t+2$ and $\dim((A_1+A_2)\cap W)\leq \ell-2n+2t+2.$ Then there exists a linearly independent subset  $\{\gamma_{t+2}, \ldots, \gamma_n, \gamma^\prime_{t+2}, \ldots, \gamma^\prime_n\}$ in $ W$ such that $
(A_1+A_2)\cap\langle\gamma_{t+2},\ldots,\gamma_n,\gamma^\prime_{t+2},\ldots,\gamma^\prime_n \rangle=\{0\}.$
Assume that
\begin{align*}
	A_1&=\langle \beta_1,\ldots,\beta_r,\beta_{r+1},\ldots,\beta_{t+1},\beta_{t+2},\ldots,\beta_{n}\rangle,\\
	A_2&=\langle \beta_1,\ldots,\beta_r,\beta^\prime_{r+1},\ldots,\beta^\prime_{t+1},\beta^\prime_{t+2},\ldots,\beta^\prime_{n}\rangle,\\
	F_1&=\langle \beta_1,\ldots,\beta_r,\beta_{r+1},\ldots,\beta_{t+1},\gamma_{t+2}+\beta_{t+2},\ldots,\gamma_{n}+\beta_{n}\rangle,\\
	F_2&=\langle \beta_1,\ldots,\beta_r,\beta^\prime_{r+1},\ldots,\beta^\prime_{t+1},\gamma^\prime_{t+2}+\beta^\prime_{t+2},\ldots,\gamma^\prime_{n}
	+\beta^\prime_{n}\rangle.
\end{align*}

Firstly, we prove that the type of $F_1$ and $F_2$ is $(n,0)$. If there exist $k_1,\ldots,k_n$ such that
\begin{align*}
\sum_{i=1}^{t+1}k_i\beta_i+\sum_{j=t+2}^{n}k_j(\gamma_j+\beta_j)=0,
\end{align*}
then we have
\begin{align*}
\sum_{i=1}^{n}k_i\beta_i+\sum_{j=t+2}^{n}k_j\gamma_j=0.
\end{align*}
It follows from $\beta_1,\ldots,\beta_n,\gamma_{t+2},\ldots,\gamma_n$ are linearly independent that $k_1=\cdots=k_n=0$. If there exist $\ell_1,\ldots,\ell_n$ such that
\begin{align*}
\sum_{i=1}^{t+1}\ell_i\beta_i+\sum_{j=t+2}^{n}\ell_j(\gamma_j+\beta_j)=\alpha\in W,
\end{align*}
then we have
\begin{align*}
\sum_{i=1}^{n}\ell_i\beta_i=\alpha-\sum_{j=t+2}^{n}\ell_j\gamma_j\in W.
\end{align*}
It follows from $A_1$ is of type $(n,0)$ that $\ell_1=\cdots=\ell_n=0$. Hence, $F_1$ is of type $(n,0)$. Similarly, we have $F_2$ is of type $(n,0)$. Since $T_1\subseteq F_1$ and $T_2\subseteq F_2$, we obtain $F_1,F_2\in\mathcal{F}$.

Next, we will prove that $\dim(F_1\cap F_2)=\dim(T_1\cap T_2)=r$. Suppose that $\dim(F_1\cap F_2)>r$. Then there exist $h_1,\ldots,h_n,h^\prime_1,\ldots,h^\prime_n$ such that
\begin{align*}
\sum_{i=1}^{t+1}h_i\beta_i+\sum_{j=t+2}^{n}h_j(\gamma_j+\beta_j)=\sum_{i=1}^{r}h_i^\prime\beta_i+\sum_{s=r+1}^{t+1}h^\prime_s\beta^\prime_s+
\sum_{j=t+2}^{n}h^\prime_j(\gamma^\prime_j+\beta^\prime_j),
\end{align*}
where $h_i\neq 0$ and $h^\prime_j\neq0$ for some $i,j\in\{r+1,r+2,\ldots,n\}$. Then
\begin{align*}
\sum_{i=1}^{r}(h_i-h^\prime_i)\beta_i+\sum_{j=r+1}^{n}h_j\beta_j-\sum_{s=r+1}^{n}h^\prime_s\beta^\prime_s+
\sum_{u=t+2}^{n}h_u\gamma_u-\sum_{v=t+2}^{n}h^\prime_v\gamma^\prime_v=0,
\end{align*}
Since $(A_1+A_2)\cap\langle\gamma_{t+2},\ldots,\gamma_n,\gamma^\prime_{t+2},\ldots,\gamma^\prime_n \rangle=\{0\}$ and $\gamma_{t+2},\ldots,\gamma_n,\gamma^\prime_{t+2},\ldots,\gamma^\prime_n$ are linearly independent, we have $h_{t+2}=\cdots=h_n=h^\prime_{t+2}=\cdots=h^\prime_n=0.$ Therefore,
\begin{align*}
\sum_{i=1}^{r}(h_i-h^\prime_i)\beta_i+\sum_{j=r+1}^{t+1}h_j\beta_j-\sum_{s=r+1}^{t+1}h^\prime_s\beta^\prime_s=0.
\end{align*}
It follows from
$\beta_1,\ldots,\beta_{t+1},\beta^\prime_{r+1},\ldots,\beta^\prime_{t+1}$ are linearly independent that $h_{r+1}=\cdots=h_{t+1}=h^\prime_{r+1}=\cdots=h^\prime_{t+1}=0$.  That is impossible since there exist some $i,j\in\{r+1,r+2,\ldots,n\}$ such that $h_i\neq 0$ and $h^\prime_j\neq0$. $\qed$
\begin{lemma}\label{b-lem3}
	Let $n,\ \ell,\ t,\ \mathcal{F}$, $\mathcal{T}$ and $M$ be as in Hypothesis~\ref{b-hyp1}. Suppose that $\mathcal{T}$ is a $t$-intersecting family, and the type of $M$ is $(e,k)$. Then the following hold.
	\begin{itemize}
		\item[{\rm(i)}] If $\tau_t(\mathcal{T})=t,$ then there exist $(t,0)$-subspace $X$ of $M$ such that $\mathcal{T}=\{T\in{M\brack t+1,0}\mid X\subseteq T\}$. Moreover, for each $F\in\mathcal{F}\setminus\mathcal{F}_X$, the type of $X+F$ is $(n+1,1)$, and $M\subseteq X+F;$ $t+1\leq e\leq n+1$ and $k\in\{0,1\}$.
		\item[{\rm(ii)}] If $\tau_t(\mathcal{T})=t+1,$ then $(e,k)$ is $(t+2,0)$ or $(t+2,1)$, and $\mathcal{T}={M\brack t+1,0}$ .
	\end{itemize}
\end{lemma}
\proof (i)\quad If $\tau_t(\mathcal{T})=t,$ then there exists a $(t,0)$-subspace $X$ such that $X\subseteq\cap_{T\in\mathcal{T}}T.$ Since $\tau_t(\mathcal{F})=t+1$, we have $\mathcal{F}\setminus\mathcal{F}_X\neq\emptyset.$ Let $F^\prime$ be an element in $\mathcal{F}\setminus\mathcal{F}_X.$ Observe that
$\dim(X\cap F^\prime)\leq t-1.$ For each $T\in\mathcal{T},$ since $X\subseteq T$ and $\dim(T\cap F^\prime)\geq t,$ we have $\dim(X\cap F^\prime)=t-1$ and $\dim(T\cap (X+F^\prime))\geq t+1,$ which implies that $X+F^\prime$ is of type $(n+1,1)$ and $T\subseteq X+F^\prime.$ Hence $M\subseteq X+F^\prime$. It is clear that $\mathcal{T}\subseteq\{T\in{M\brack t+1,0}\mid X\subseteq T\}$. Let $T^\prime$ be a $(t+1,0)$-subspace of $V$ with $X\subseteq T^\prime\subseteq M.$ For each $F\in\mathcal{F}_X$, then $\dim(T^\prime\cap F)\geq t$. For each $F\in\mathcal{F}\setminus\mathcal{F}_X$,  by the discussion above, then $T^\prime\subseteq M\subseteq X+F,$ implying that $\dim(T^\prime\cap F)\geq t$. Hence, $T^\prime\in\mathcal{T}$, and (i) holds.
\
\medskip
\

(ii)\quad Since $\tau_t(\mathcal{T})=t+1$, by the structure of maximal $t$-intersecting families of ${V\brack t+1,0}$, we have $|\mathcal{T}|\geq 3$, $M$ is of type $(t+2,0)$ or $(t+2,1)$, and $\mathcal{T}\subseteq{M\brack t+1,0}$. Note that there exist three distinct subspaces $T_1,T_2,T_3\in\mathcal{T}$ such that $T_1\cap T_2$, $T_1\cap T_3$ and $T_2\cap T_3$ are pairwise distinct.
 For each $F\in \mathcal{F},$ if $F\cap T_1=F\cap T_2=F\cap T_3,$ then $F\cap T_1\subseteq T_i$ for each $i\in\{1,2,3\}.$ That is impossible because $T_1\cap T_2,$ $T_1\cap T_3$ and $T_2\cap T_3$ are pairwise distinct and $\dim(F\cap T_1)\geq t$. Hence there exist $T_i,T_j\in\{T_1,T_2,T_3\}$ such that $F\cap T_i\neq F\cap T_j,$ implying that $\dim(F\cap M)\geq t+1.$ Then for each $T^\prime\in{M\brack t,0}$, we have $\dim(F\cap T^\prime)\geq t.$ Therefore $\mathcal{T}={M\brack t+1,0}$ as desired. $\qed$
\begin{lemma}\label{b-prop3}
	Let $n,\ \ell,\ t,\ \mathcal{F}$, $\mathcal{T}$ and $M$ be as in Hypothesis~\ref{b-hyp1}. Suppose $\mathcal{T}$ is a $t$-intersecting family with $\tau_t(\mathcal{T})=t$,  $X$ is a $(t,0)$-subspace of  $\cap_{T\in\mathcal{T}}T,$ and the type of $M$ is $(n+1,1)$. Then $\mathcal{F}=\mathcal{H}_1(X,M).$
\end{lemma}
\proof By Lemma~\ref{b-lem3} (i), for each $F\in\mathcal{F}\setminus\mathcal{F}_X$,  we have $M=F+X,$ implying that $F\in{M\brack n,0}.$  Let $\mathcal{A}^\prime=\{F\in {V\brack n,0}\mid X\subseteq F,\ \dim(F\cap M)\geq t+1\}$ and $F^\prime$ a fixed element in $\mathcal{F}\setminus\mathcal{F}_X$. For each $F\in\mathcal{F}_X,$ since $\dim(F\cap F^\prime)\geq t,$ $\dim(F^\prime\cap X)=t-1$ and $M=F^\prime+X,$ we have $\dim(F\cap M)\geq t+1,$ which implies that $\mathcal{F}_X\subseteq \mathcal{A}^\prime.$ Note that $\mathcal{A}^\prime\cup {M\brack n,0}$ is a $t$-intersecting family. By the maximality of $\mathcal{F},$ we have $\mathcal{F}=\mathcal{A}^\prime\cup {M\brack n,0}.$ $\qed$
\begin{lemma}\label{b-prop3-1}
	Let $n,\ \ell,\ t,\ \mathcal{F}$, $\mathcal{T}$ and $M$ be as in Hypothesis~\ref{b-hyp1}. Suppose that $\mathcal{T}$ is a $t$-intersecting family with $\tau_t(\mathcal{T})=t$, $X$ is a $(t,0)$-subspace of $\cap_{T\in\mathcal{T}}T,$ and the type of $M$ is $(n,k)$, where $k\in\{0,1\}$. Assume that $C=M+\sum\limits_{F\in\mathcal{F}\setminus \mathcal{F}_X}F$ and the type of $C$ is $(c,c-n)$. Then $n+2\leq c\leq 2n-t$ or $c=n+\ell.$ Moreover, the following hold:
	\begin{itemize}
		\item[{\rm(i)}]If $n+2\leq c\leq 2n-t,$ then $\mathcal{F}=\mathcal{A}(X,M,C)\cup\mathcal{B}(X,M,C)\cup\mathcal{C}(X,M,C),$
		\item[{\rm(ii)}] If $c=n+\ell$, then  $\mathcal{F}=\mathcal{A}(X,M,V)\cup\mathcal{C}(X,M,V),$
	\end{itemize}
where $\mathcal{A}(X,M,C),\mathcal{B}(X,M,C),\mathcal{C}(X,M,C)$ are as in Family~II.
\end{lemma}
\proof By Lemma~\ref{b-lem3} (i), for each $F\in\mathcal{F}\setminus\mathcal{F}_X,$ we have $\dim(F\cap X)=t-1$ and $M\subseteq X+F$. It follows from $X\subseteq M$ that $\dim(F\cap M)=n-1$. Similar to the proof of Lemma~3.5 in \cite{Cao-vec}, we can choose $F_1,F_2,\ldots, F_{c-n}\in\mathcal{F}\setminus\mathcal{F}_X$ such that $F_i\cap(M+\sum\limits_{j=1}^{i-1}F_j)=F_i\cap M$ for each $i\in\{1,2,\ldots,c-n\}.$ If there exists $F^\prime\in\mathcal{F}$ such that $F^\prime\cap M=X$, then for every $i\in\{1,2,\ldots,c-n\}$, there exists $y_i\in F_i\setminus M$ such that $y_i\in F^\prime$ due to $\dim(F^\prime \cap F_i)\geq t$ and $\dim(F^\prime \cap F_i\cap M)=t-1.$ Let $x_1,\ldots,x_t$ be a basis of $X$. By the choices of $F_1,F_2,\ldots,F_{c-n},$ it is routine to prove that $x_1,\ldots,x_t,y_1\ldots,y_{c-n}$ are linearly independent in $F^\prime$ since $F_i\cap(M+\sum\limits_{j=1}^{i-1}F_j)=F_i\cap M$ for every $i\in\{1,2,\ldots,c-n\}.$

 Suppose $c=n+1.$ Observe that $C=F^\prime+X$ for each $F^\prime\in\mathcal{F}\setminus\mathcal{F}_X$. Then for each $T\in{C\brack t+1,0}$ with $X\subseteq T$, we have $\dim(T\cap F)\geq t$ for all $F\in\mathcal{F}$, which implies that $T\subseteq M$, a contradiction.

Suppose $c\geq 2n-t+1.$ If there exists $F^\prime\in\mathcal{F}$ such that $F^\prime\cap M=X$, by the discussion above, we can obtain $c-n+t$ vectors in $F^\prime$ which are linearly independent. That is impossible. Hence, we have $\dim(F^{\prime\prime}\cap M)\geq t+1$ for all $F^{\prime\prime}\in\mathcal{F}_X$. By the maximality of $\mathcal{F}$, it is easy to see that each $(n,0)$-subspace $F^{\prime\prime\prime}$ of $V$ satisfying $\dim(F^{\prime\prime\prime}\cap X)=t-1$ and $\dim(F^{\prime\prime\prime}\cap M)=n-1$ is in $\mathcal{F}$. Let $\alpha_1,\alpha_2,\ldots,\alpha_\ell$ be a basis of $W$, and $\alpha_0=0$.

Assume that $M$ is of type $(n,0)$. Let $H$ be an $(n-1,0)$-subspace  of $V$ such that $X\nsubseteq H\subseteq M$, and $\alpha$ a vector in $M$ with $\alpha\notin X$ and $\alpha\notin H$. Suppose $H_i=H+\langle \alpha_i+\alpha\rangle$ for each $i\in\{1,2,\ldots,\ell\}$. Note that $\alpha_i+\alpha\notin M$ and $\alpha_i+\alpha\notin W$. Then $\dim(H_i)=n,$ $\dim(H_i\cap M)=n-1$ and $\dim(H_i\cap X)=t-1$. If there exists $\gamma\in H_i\cap W$, then there exist $\gamma^\prime\in H$ and $k_i\in\mathbb{F}_q$ such that $\gamma=\gamma^\prime+k_i(\alpha_i+\alpha)\in W$. Then $\gamma^\prime+k_i\alpha\in W,$ implying that $\gamma^\prime=0$ and $k_i=0$. Hence, the type of $H_i$ is $(n,0)$. Then $V=M+\sum_{i=1}^\ell H_i\subseteq C$ and $c=n+\ell$.

Assume that $M$ is of type $(n,1)$. Observe that the number of $(n-1,0)$-subspaces $H$ of $V$ satisfying $X\nsubseteq H\subseteq M$ is $N^\prime(0,0;n-1,0;n,n-1)-N^{\prime}(t,0;n-1,0;n,n-1)=q^{n-1}-q^{n-1-t}>0.$ Suppose $H^\prime$ is one of that $(n-1,0)$-subspaces, and $\alpha$ is a vector in $V\setminus(W+M)$. Assume that $H_i=H^\prime+\langle \alpha_i+\alpha\rangle$ for each $i\in\{0,1,\ldots,\ell\}$. For every $i\in\{1,2,\ldots,\ell\}$, we have that if $\alpha_i+\alpha\in H^\prime$, then $\alpha_i\in H^\prime+\langle \alpha\rangle$ and $\dim(W\cap (H^\prime+\langle\alpha\rangle))\geq1$, which contradict with $V=W+H^\prime+\langle\alpha\rangle$. If $\beta\in H_i\cap W,$ then there exist $\gamma^\prime\in H^\prime$ and $k_i\in\mathbb{F}_q$ such that $\gamma^\prime+k_i(\alpha_i+\alpha)=\beta\in W,$ implying that $\gamma^\prime+k_i\alpha=\beta-k_i\alpha_i\in W$. Hence, $\gamma^\prime=0$, $k_i=0$ and $\beta=0$. Therefore $H_i$ is of type $(n,0)$. Then $V=M+\sum_{i=0}^\ell H_i\subseteq C$ and $c=n+\ell.$
\
\medskip
\

(i)\quad Suppose $n+2\leq c\leq 2n-t$.  Since  $\dim(F\cap X)=t-1$ and $\dim(F\cap M)=n-1$ for each $F\in\mathcal{F}\setminus\mathcal{F}_X,$ we have $\mathcal{F}\setminus\mathcal{F}_X\subseteq \mathcal{C}(X,M,C).$  For each $F^\prime\in\mathcal{F}_X,$ if $\dim(F^\prime\cap M)\geq t+1$, then $F^\prime\in\mathcal{A}(X,M,C);$ if $F^\prime\cap M=X,$ then $\dim(F^\prime\cap C)=c-n+t$ by the discussion above, which implies that $F^\prime\in\mathcal{B}(X,M,C).$ Thus, $\mathcal{F}\subseteq \mathcal{A}(X,M,C)\cup\mathcal{B}(X,M,C)\cup\mathcal{C}(X,M,C).$ It is routine to check that $\mathcal{A}(X,M,C)\cup\mathcal{B}(X,M,C)\cup\mathcal{C}(X,M,C)$ is a $t$-intersecting family. By the maximality of $\mathcal{F}$, we get $\mathcal{F}=\mathcal{A}(X,M,C)\cup\mathcal{B}(X,M,C)\cup\mathcal{C}(X,M,C)$ as required.
\
\medskip
\

(ii)\quad  In this case, by the discussion in (i) and  maximality of $\mathcal{F}$ again, observe that $\mathcal{F}=\mathcal{A}(X,M,C)\cup\mathcal{C}(X,M,C).$  $\qed$
\begin{lemma}\label{b-prop4}
	Let $n,\ \ell,\ t,\ \mathcal{F}$, $\mathcal{T}$  and $M$ be as in Hypothesis~\ref{b-hyp1}. If $\mathcal{T}$ is a $t$-intersecting family with $\tau_t(\mathcal{T})=t+1$, then $\mathcal{F}=\mathcal{H}_3(M).$
\end{lemma}
\proof Suppose that the type of $M$ is $(t+2,k)$, where $k\in\{0,1\}.$ Since $\mathcal{T}={M\brack t+1,0}$, we have $\dim(F\cap M)\geq t$ for all $F\in\mathcal{F}.$ If there exists $F^\prime\in\mathcal{F}$ such that $\dim(F^\prime\cap M)=t,$ from
\begin{align*}
\left|\left\{T\in{M\brack t+1,0}\mid F^\prime\cap M\nsubseteq T\right\}\right|=\left|{M\brack t+1,0}\right|-N^\prime(t,0;t+1,0;t+2,t+2-k)>0,
\end{align*}
then there exists a $T^\prime\in\mathcal{T}$ such that $\dim(F^\prime\cap T^\prime)=t-1,$ a contradiction. Hence, $\mathcal{F}\subseteq \mathcal{H}_3(M).$ Since $\mathcal{F}$ is maximal and $\mathcal{H}_3(M)$ is $t$-intersecting, we have $\mathcal{F}=\mathcal{H}_3(M)$ as desired. $\qed$
\begin{lemma}\label{b-prop5}
	Let $n,\ \ell,\ t,\ \mathcal{F}$ and $\mathcal{T}$ be as in Hypothesis~\ref{b-hyp1}. If $\mathcal{T}$ is not $t$-intersecting, then
	\begin{align*}
	|\mathcal{F}|\leq (q+1)^2\cdot q^{\ell(n-t-1)}+(q^2+2q+2)\theta_{t-1}\cdot q^{\ell(n-t-2)+2}+\theta_{t-1}\theta_{t-2}\cdot q^{\ell(n-t-3)+5}.
\end{align*}
\end{lemma}
\proof Since $\mathcal{T}$ is not $t$-intersecting, there exist $T_1$ and $T_2$ in $\mathcal{T}$ such that $\dim(T_1\cap T_2)\leq t-1$. Suppose $\dim(T_1\cap T_2)=d$. For each $F\in\mathcal{F},$ since $\dim(F\cap T_i)\geq t$ for every $i\in\{1,2\},$ there exist $E_1\in{T_1\brack t}$ and $E_2\in{T_2\brack t}$ such that $E_1\subseteq F$ and $E_2\subseteq F$, implying that
\begin{align*}
\mathcal{F}=\bigcup_{E_1\in{T_1\brack t},\ E_2\in{T_2\brack t}} \mathcal{F}_{E_1+E_2}.
\end{align*}

Let $E_1\in{T_1\brack t}$ and $E_2\in{T_2\brack t}$ be two subspaces. Suppose the type of $E_1+E_2$ is $(e,k)$. Since $E_1\cap E_2\subseteq T_1\cap T_2$ and $E_1+E_2\subseteq T_1+T_2$, we have $\min\{2t+2-d,2t\}\geq e=\dim(E_1+E_2)=2t-\dim(E_1\cap E_2)\geq 2t-d$ and
\begin{align*}
|\mathcal{F}_{E_1+E_2}|\leq N^\prime(e,k;n,0;n+\ell,n)\leq q^{\ell(n-e)}.
\end{align*}
For every $e\in\{2t-d,2t-d+1,\ldots,\min\{2t+2-d,2t\}\}$, let
\begin{align*}
\mathcal{E}_e(T_1,T_2)=\left\{(E_1,E_2)\in {T_1\brack t}\times{T_2\brack t}\mid \dim(E_1+E_2)=e\right\}.
\end{align*}
Then
\begin{align}\label{b-equ7}
	|\mathcal{F}|\leq \sum_{e=2t-d}^{\min\{2t+2-d,2t\}}\left|\mathcal{E}_e(T_1,T_2)\right|\cdot q^{\ell(n-e)}.
\end{align}

If $e=2t-d,$ then $E_1\cap E_2=T_1\cap T_2$ for each $(E_1, E_2)\in \mathcal{E}_{2t-d}(T_1, T_2)$, and
\begin{align*}
\mathcal{E}_{2t-d}(T_1,T_2)=\left\{(E_1,E_2)\in {T_1\brack t}\times{T_2\brack t}\mid T_1\cap T_2\subseteq E_1,\ T_1\cap T_2\subseteq E_2\right\},
\end{align*}
implying that $|\mathcal{E}_{2t-d}(T_1,T_2)|=(\theta_{t+1-d})^2$.

If $e=2t-d+1,$ then $\dim(E_1\cap E_2)=\dim(T_1\cap T_2)-1$ for each $(E_1, E_2)\in \mathcal{E}_{2t-d+1}(T_1, T_2)$, and
\begin{align*}
	\mathcal{E}_{2t-d+1}(T_1,T_2)=\bigcup_{D\in{T_1\cap T_2\brack d-1}}\Big{(}&\left\{(E_1,E_2)\in {T_1\brack t}\times{T_2\brack t}\mid
	D\subseteq E_1,\ D\subseteq E_2\right\}\\
	&\setminus\left\{(E_1,E_2)\in {T_1\brack t}\times{T_2\brack t}\mid
	T_1\cap T_2\subseteq E_1,\ T_1\cap T_2\subseteq E_2\right\}\Big{)},
\end{align*}
implying that
\begin{align*}
|\mathcal{E}_{2t-d+1}(T_1,T_2)|=\theta_d\left((\theta_{t+2-d})^2-(\theta_{t+1-d})^2\right).
\end{align*}

If $e=2t-d+2,$ then $|\mathcal{E}_{2t-d+2}(T_1,T_2)|=(\theta_{t+1})^2-|\mathcal{E}_{2t-d+1}(T_1,T_2)|-|\mathcal{E}_{2t-d}(T_1,T_2)|$ from $\min\{2t+2-d,2t\}\geq e\geq 2t-d$.

Therefore, by (\ref{b-equ7}), we have
\begin{align*}
	|\mathcal{F}|\leq& (\theta_{t+1-d})^2\cdot q^{\ell(n-2t+d)}+ \theta_d\left((\theta_{t+2-d})^2-(\theta_{t+1-d})^2\right)\cdot
 q^{\ell(n-2t+d-1)}\\ &+\left((\theta_{t+1})^2- \theta_d(\theta_{t+2-d})^2+(\theta_d-1)(\theta_{t+1-d})^2\right)\cdot q^{\ell(n-2t+d-2)}.
\end{align*}
It is routine to verify that the three functions $(\theta_{t+1-d})^2\cdot q^{\ell(n-2t+d)}$, $\theta_d\left((\theta_{t+2-d})^2-(\theta_{t+1-d})^2\right)\cdot q^{\ell(n-2t+d-1)}$ and $(\theta_{t+1})^2- \theta_d(\theta_{t+2-d})^2+(\theta_d-1)(\theta_{t+1-d})^2$ are all increasing as $d\in\{0,1,\ldots,t-1\}$ increases. Hence, the required result follows by setting $d=t-1$.  $\qed$
\begin{lemma}\label{b-prop1}
	Let $n,\ \ell,\ t,\ \mathcal{F},\ \mathcal{T}$ and $M$ be as in Hypothesis~\ref{b-hyp1}. Suppose that $\mathcal{T}$ is a $t$-intersecting family with $\tau_t(\mathcal{T})=t.$ Then the following hold.
	\begin{itemize}
		\item[{\rm (i)}] If $\dim(M)=t+1$, then
		\begin{align*}
			|\mathcal{F}|\leq q^{\ell(n-t-1)}+\left(q\theta_{n-t}-1\right)q^{\ell(n-t-2)+2}\theta_{t+1}\theta_{n-t-1}.
    	\end{align*}	
		\item[{\rm (ii)}] If $\dim(M)=t+2\leq n-1,$ then
		\begin{align*}
			|\mathcal{F}|\leq q^{\ell(n-t-1)}(q+1)+q^{\ell(n-t-2)+2}(q^{n-t+1}+q^{t+2}-2q)\theta_{n-t-1}.
		\end{align*}
		\item[{\rm (iii)}] If $t+3\leq\dim(M)\leq n-1,$ then
		\begin{align*}
			|\mathcal{F}|\leq q^{\ell(n-t-1)}\theta_{n-t-1}+q^{\ell(n-t-2)+3}\theta_{n-t}\theta_{n-t-1}+q^{\ell(n-t-2)+t+3}.
		\end{align*}
	\end{itemize}
\end{lemma}
\proof (i) Observe that $M$ is the unique $(t+1,0)$-subspace of $V$ in $\mathcal{T}$. Since $\dim(M\cap F)\geq t$ for all $F\in\mathcal{F}$, we have
\begin{eqnarray}\label{b-equ1-1}
	\mathcal{F}=\mathcal{F}_M\cup\left(\bigcup_{S\in{M\brack t}}(\mathcal{F}_S\setminus\mathcal{F}_M)\right).
\end{eqnarray}
Next, we give an upper bound of $|\mathcal{F}_S\setminus\mathcal{F}_M|$ for each $S\in{M\brack t}.$

Let $S$ be a $t$-subspace of $M$. Since $\tau_t(\mathcal{F})=t+1,$ there exists an $F^\prime\in\mathcal{F}\setminus \mathcal{F}_S$ satisfying $\dim(S\cap F^\prime)=t-1$ due to $\dim(F^\prime \cap M)\geq t.$ Then we have  $M=(F^\prime \cap M)+S$ and $M\subseteq F^\prime +S.$ Set
\begin{align*}
\mathcal{W}=\left\{H\in{F^\prime+S\brack t+1,0}\mid S\subseteq H,\ H\neq M\right\}
\end{align*}
For each $F\in \mathcal{F}_S\setminus\mathcal{F}_M$, notice that  $\dim(F\cap(F^\prime+S))=\dim(F)+\dim(F^\prime+S)-\dim(F+F^\prime+S)\geq t+1.$ Hence there exists an $H\in\mathcal{W}$ such that $H\subseteq F.$ Therefore, we have
\begin{eqnarray}\label{b-equ2-1}
	\mathcal{F}_S\setminus\mathcal{F}_T=\bigcup_{H\in\mathcal{W}} \mathcal{F}_H.
\end{eqnarray}

For each $H\in\mathcal{W},$ since $\mathcal{T}=\{M\}$, there exists $F^{\prime\prime}$ such that $\dim(H\cap F^{\prime\prime})<t,$ which implies that $\dim(H\cap F^{\prime\prime})=t-1$ due to $\dim(H\cap M)=\dim(S)=t$ and $\dim(M\cap F^{\prime\prime})\geq t.$ From Lemma~\ref{b-lem8}, we have $|\mathcal{F}_{H}|\leq q^{\ell(n-t-2)+2}\theta_{n-t-1}.$ By Lemma~\ref{lem5}, we obtain $|\mathcal{F}_M|\leq N^\prime(t+1,0;n,0;n+\ell,n)=q^{\ell(n-t-1)},$ and
$\left|\mathcal{W}\right|=N^\prime(t,0;t+1,0;n+1,n)-1=q\theta_{n-t}-1$
since $S+F^\prime$ is of type $(n+1,1)$. Therefore, from (\ref{b-equ1-1}) and (\ref{b-equ2-1}), we obtain
\begin{eqnarray*}
	|\mathcal{F}|\leq q^{\ell(n-t-1)}+\theta_{t+1}\cdot \left(q\theta_{n-t}-1\right)\cdot q^{\ell(n-t-2)+2}\theta_{n-t-1},
\end{eqnarray*}
as required.
\
\medskip
\

(ii) and (iii)\quad Suppose $X=\cap_{T\in\mathcal{T}}T$, and the type of $M$ is $(e,k)$ with $t+2\leq e\leq n-1$ and $k\in\{0,1\}$. Since $|\mathcal{T}|\geq 2$ and $\tau_{t}(\mathcal{T})=t$, we have $\dim(X)=t.$ In the following, we will obtain the upper bound of $|\mathcal{F}|$ by giving the upper bound of $|\mathcal{F}_X|$ and $|\mathcal{F}\setminus\mathcal{F}_X|$ respectively. Since $\tau_t(\mathcal{F})=t+1,$ we have $\dim(F\cap X)\geq t-1$ for all $F\in\mathcal{F}$, and there exists $F^\prime\in\mathcal{F}$ such that $\dim(X\cap F^\prime)=t-1.$ From Lemma~\ref{b-lem3} (i), we have $X\subseteq M\subseteq X+ F^\prime.$ Set
\begin{align*}
\mathcal{W}_1=\left\{H\in{M\brack t+1,0}\mid X\subseteq H\right\}\quad \mbox{and}\quad \mathcal{W}_2=\left\{H\in{X+F^\prime\brack t+1,0}\mid X\subseteq H\right\}\setminus\mathcal{W}_1.
\end{align*}

For each $F\in\mathcal{F}_X,$ notice that $\dim(F\cap(X+ F^\prime))\geq t+1$ due to $X\subseteq F$ and $\dim(F\cap F^\prime)\geq t.$ Then we have
\begin{eqnarray*}
	\mathcal{F}_X=\left(\bigcup_ {H_1\in\mathcal{W}_1}\mathcal{F}_{H_1}\right)\cup\left(\bigcup_{H_2\in\mathcal{W}_2}\mathcal{F}_{H_2}\right).
\end{eqnarray*}
From Lemma~\ref{lem5}, we have $|\mathcal{F}_{H_1}|\leq N^\prime(t+1,0;n,0;n+\ell,n)=q^{\ell(n-t-1)}$ for each $H_1\in\mathcal{W}_1$, and  $|\mathcal{W}_1|=N^\prime(t,0;t+1,0;e,e-k)=q^k\theta_{e-k-t},$  implying that $|\bigcup_{H_1\in\mathcal{W}_1}\mathcal{F}_{H_1}|\leq q^{k+\ell(n-t-1)}\theta_{e-k-t}.$ For each $H_2\in\mathcal{W}_2,$ since $H_2\notin \mathcal{T},$ there exists $F^{\prime\prime}\in\mathcal{F}$ such that $\dim(H_2\cap F^{\prime\prime})<t.$ It follows from $\dim(F^{\prime\prime}\cap X)\geq t-1$ that $\dim(H_2\cap F^{\prime\prime})=t-1$. Then $|\mathcal{F}_{H_2}|\leq q^{\ell(n-t-2)+2}\theta_{n-t-1}$ by Lemma~\ref{b-lem8}. Notice that
\begin{align*}
	\left|\mathcal{W}_2\right|=N^\prime(t,0;t+1,0;n+1,n)-N^\prime(t,0;t+1,0;e,e-k)=q\theta_{n-t}-q^k\theta_{e-k-t}
\end{align*}
by Lemma~\ref{lem5}. Therefore, we have
\begin{eqnarray}\label{b-equ11}
	|\mathcal{F}_X|\leq q^{k+\ell(n-t-1)}\theta_{e-k-t}+\left(q\theta_{n-t}-q^k\theta_{e-k-t}\right)\cdot q^{\ell(n-t-2)+2}\theta_{n-t-1}.
\end{eqnarray}

For each $F\in\mathcal{F}\setminus\mathcal{F}_X$, by Lemma~\ref{b-lem3} (i), we have $\dim(F\cap X)=t-1$ and $M\subseteq F+X,$ which implies that $\dim(M\cap F)=e-1$. Set
\begin{align*}
\mathcal{W}_3=\left\{D\in{M\brack e-1,0}\mid X\nsubseteq D\right\}.
\end{align*}
Hence,
\begin{align*}
	\mathcal{F}\setminus\mathcal{F}_X\subseteq&\left\{F\in{V\brack n,0}\mid\dim(F\cap M)=e-1,\ X\nsubseteq F\right\}\\
	=&\bigcup_{D\in\mathcal{W}_3}\left(\left\{F\in{V\brack n,0}\mid D\subseteq F\right\}\setminus \left\{F\in{V\brack n,0}\mid M\subseteq F\right\}\right).
\end{align*}
Observe that the number of $(n,0)$-subspaces $F$ of $V$ with $D\subseteq F$ is $N^\prime(e-1,0;n,0;n+\ell,n)=q^{\ell(n-e+1)}$ for each $D\in{M\brack e-1,0},$ and the number of $(n,0)$-subspaces $F$ of $V$ with $M\subseteq F$ is $N^\prime(e,k;n,0;n+\ell,n)=q^{\ell(n-e+k)}{\ell-k\brack -k}$. Since  $|\mathcal{W}_3|=N^\prime(0,0;e-1,0;e,e-k)-N^\prime(t,0;e-1,0;e,e-k)=q^{k(e-1)}{e-k\brack e-1}-q^{k(e-1-t)}{e-k-t\brack e-1-t}$, we have
\begin{eqnarray}\label{b-equ12}
	|\mathcal{F}\setminus\mathcal{F}_X|\leq \left(q^{k(e-1)}{e-k\brack e-1}-q^{k(e-1-t)}{e-k-t\brack e-1-t}\right)\cdot\left(q^{\ell(n-e+1)}-q^{\ell(n-e+k)}{\ell-k\brack -k}\right).
\end{eqnarray}

In particular, we consider the case when $e=t+2$. By the discussion above, we have $\dim(M\cap F)=e-1=t+1$ for all $F\in\mathcal{F}\setminus\mathcal{F}_X,$ which implies that
\begin{align*}
\mathcal{F}\setminus\mathcal{F}_X\subseteq \bigcup_{L\in\mathcal{W}_3} \mathcal{F}_L.
\end{align*}
For each $L\in\mathcal{W}_3$, since $L\notin\mathcal{T}$ and $\dim(F\cap M)\geq t$ for all $F\in\mathcal{F}$, there exists $F^\prime\in\mathcal{F}$ such that $\dim(F^\prime\cap L)=t-1.$ Then $|\mathcal{F}_L|\leq q^{\ell(n-t-2)+2}\theta_{n-t-1}$ by Lemma~\ref{b-lem8}, and
\begin{eqnarray}\label{b-equ13}
	|\mathcal{F}\setminus\mathcal{F}_X|\leq \left(q^{k(t+1)}{t+2-k\brack 1-k}-q^k\theta_{2-k}\right)\cdot q^{\ell(n-t-2)+2}\theta_{n-t-1}.
\end{eqnarray}

If $\dim(M)=t+2$, by  (\ref{b-equ11}) and (\ref{b-equ13}), then
\begin{align*}
	|\mathcal{F}|\leq& q^{k+\ell(n-t-1)}\theta_{2-k}+\left(q\theta_{n-t}-2q^k\theta_{2-k}+q^{k(t+1)}{t+2-k\brack 1-k}\right)q^{\ell(n-t+2)-2}\theta_{n-t-1}\\
	<&q^{\ell(n-t-1)}(q+1)+(q^{n-t+1}-2q+q^{t+2})q^{\ell(n-t+2)-2}\theta_{n-t-1},
\end{align*}
and (ii) holds.

Suppose $\dim(M)\geq t+3$. Since $q^{k(e-1)}{e-k\brack e-1}-q^{k(e-1-t)}{e-k-t\brack e-1-t}<q^e$ for $k\in\{0,1\}$, by (\ref{b-equ11}) and (\ref{b-equ12}), we have
\begin{align*}
	|\mathcal{F}|<& q^{\ell(n-t-1)}\theta_{e-t}+q^{\ell(n-t-2)+3}\theta_{n-t}\theta_{n-t-1}+q^{\ell(n-e+1)+e} \\
	\leq&q^{\ell(n-t-1)}\theta_{n-1-t}+q^{\ell(n-t-2)+3}\theta_{n-t}\theta_{n-t-1}+q^{\ell(n-t-2)+t+3},
\end{align*}
and (iii) holds.    $\qed$
\section{The proof of Theorem~\ref{b-main1}}
In this section, we will prove Theorem~\ref{b-main1}. Firstly, we give some inequalities for the sizes of the families in Families~I, II and III. Set
\begin{align*}
f^\prime(n,\ell,t)=q^{\ell(n-t-1)+1}\theta_{n-t-1}-q^{\ell(n-t-2)+3}{n-t-1\brack 2}.
\end{align*}

\begin{lemma}\label{b-lem6}
Suppose $4\leq n+1\leq \ell$ and $1\leq t\leq n-2.$  Let $\mathcal{H}_2(X,M,C)$ be a family constructed in Family II. Then $|\mathcal{H}_2(X,M,C)|\geq f^\prime(n,\ell,t)$. Moreover,  if $1\leq t\leq k-3$ and $\dim(C)=n+1$, then
$|\mathcal{H}_2(X,M,C)|\leq q^{\ell(n-t-1)+1}\theta_{n-t}.$
\end{lemma}
\proof Let $N\in{V\brack e,r}$ with $X\subseteq N$, where $t\leq e$ and $r\in\{0,1\}$. For $i\in\{t,t+1,\ldots,n\}$, set
\begin{align*}
\mathcal{A}_i(X,N)=&\left\{F\in{V\brack n,0}\mid X\subseteq F,\ \dim(F\cap N)=i\right\},\\
\mathcal{L}_i(X,N)=&\left\{(I,F)\in{V\brack i,0}\times{V\brack n,0}\mid X\subseteq I\subseteq N,\ I\subseteq F \right\}.
\end{align*}
Counting $|\mathcal{L}_i(X,N)|$ in two ways, by Lemma \ref{lem5}, we have
\begin{align}
|\mathcal{L}_i(X,N)|=&\sum_{j=i}^n\left|\mathcal{A}_j(X,N)\right|\cdot{j-t\brack i-t}\nonumber\\
=&N^\prime(t,0;i,0;e,e-r)N^\prime(i,0;n,0;n+\ell,n)=q^{r(i-t)+\ell(n-i)}{e-r-t\brack i-t}.\label{b-equ15}
\end{align}
In particular, we have
\begin{eqnarray}
|\mathcal{L}_{t+1}(X,N)|=\sum_{j=t+1}^n\left|\mathcal{A}_j(X,N)\right|+\sum_{j=t+2}^n\left|\mathcal{A}_j(X,N)\right|\cdot\left(\theta_{j-t}-1\right)\label{b-equ14}.
\end{eqnarray}

Suppose that $\mathcal{A}(X,M,C)$ is the family constructed in Family II. Observe that $(e,r)=(n,k)$, and $\mathcal{A}(X,M,C)=\cup_{j=t+1}^n\mathcal{A}_j(X,M).$ By (\ref{b-equ15}) and (\ref{b-equ14}), we obtain
\begin{align*}
|\mathcal{L}_{t+1}(X,M)|\leq|\mathcal{A}(X,M,C)|+\sum_{j=t+2}^k\left|\mathcal{A}_j(X,M)\right|\cdot q{j-t\brack 2}
=|\mathcal{A}(X,M,C)|+q|\mathcal{L}_{t+2}(X,M)|.
\end{align*}
By (\ref{b-equ15}) again, it follows that
\begin{align*}
|\mathcal{A}(X,M,C)|\geq\left\{
\begin{array}{ll}
	q^{\ell(n-t-1)}\theta_{n-t}-q^{\ell(n-t-2)+1}{n-t\brack 2}, & \mbox{if}\ k=0,\\
q^{\ell(n-t-1)+1}\theta_{n-t-1}-q^{\ell(n-t-2)+3}{n-t-1\brack 2}, &	 \mbox{if}\ k=1.
\end{array}
\right.
\end{align*}
Observe that the lower bound given in the second line of the display above is $f^\prime(n,\ell,t)$, and
\begin{align*}
&q^{\ell(n-t-1)}\theta_{n-t}-q^{\ell(n-t-2)+1}{n-t\brack 2}-f^\prime(n,\ell,t)\\
=&q^{\ell(n-t-1)}-q^{\ell(n-t-2)+1}\theta_{n-t-1}	\geq q^{\ell(n-t-1)}-q^{\ell(n-t-2)+n-t}>0.
\end{align*}	
From the construction of $\mathcal{H}_2(X,M,C)$, we have $|\mathcal{H}_2(X,M,C)|\geq f^\prime(n,\ell,t)$ as required.

 Suppose $C$ is of type $(n+1,1)$, and set $\mathcal{A}(X,C)=\cup_{j=t+1}^n\mathcal{A}_j(X,C).$ Observe that
\begin{align*}
\mathcal{A}_{n-1}(X,C)=\bigcup_{D\in{C\brack n-1,0},\ X\subseteq D}\left\{F\in{V\brack n,0}\mid D\subseteq F,\ F\nsubseteq C\right\}.
\end{align*}
Then by Lemmas~\ref{lem5}, we have
$
|\mathcal{A}_{n-1}(X,C)|=N^\prime(t,0;n-1,0;n+1,n)\cdot (N^\prime(n-1,0;n,0;n+\ell,n)-N^\prime(n-1,0;n,0;n+1,n))
$
and
\begin{align*}
	\left|\mathcal{A}_{n-1}(X,C)\right|\cdot\left(\theta_{n-1-t}-1\right)
	=q^{n-t}(q^{\ell-1}-1)\theta_{n-t}\cdot q\theta_{n-t-2}
	>q^n-q^{n-t}.
\end{align*}

For each $F\in{C\brack n,0}$, observe that $F\in\mathcal{A}(X,C)$ if and only if $X\subseteq F$. By Lemma~\ref{lem5}, then
\begin{align}\label{n-1}
\left|{C\brack n,0}\setminus \mathcal{A}(X,C)\right|=N^\prime(0,0;n,0;n+1,n)-N^\prime(t,0;n,0;n+1,n)=q^n-q^{n-t}.
\end{align}

Since $1\leq t\leq n-3$, by (\ref{b-equ14}), we have
\begin{align*}
	|\mathcal{L}_{t+1}(X,C)|=&|\mathcal{A}(X,C)|+\sum_{j=t+2}^n\left|\mathcal{A}_j(X,C)\right|\cdot \left(\theta_{j-t}-1\right)\\
	\geq&|\mathcal{A}(X,C)|+\left|\mathcal{A}_{n-1}(X,C)\right|\cdot \left(\theta_{n-1-t}-1\right)
	\geq|\mathcal{A}(X,C)|+q^n-q^{n-t}.
\end{align*}
 By (\ref{b-equ15}), Remark~\ref{b-obs} and the structure of Family I,   $|\mathcal{H}_2(X,M,C)|\leq q^{\ell(n-t-1)+1}\theta_{n-t}$.  $\qed$
\begin{lemma}\label{b-lll2}
The size of each family constructed in Family II only depends on the type of $Z$. Set $h_3(t+2,k)=|\mathcal{H}_3(Z)|$ for $(t+2,k)$-subspace $Z$ of $V$. Then the following hold:
\begin{align}
	h_3(t+2,0)=&q^{\ell(n-t-1)}\theta_{t+2}-q^{\ell(n-t-2)+1}\theta_{t+1},\label{b-equ1}\\
	h_3(t+2,1)=& q^{\ell(n-t-1)+t+1}.\label{b-equ10}
\end{align}
Moreover, $h_3(t+2,0)>h_3(t+2,1).$
\end{lemma}
\proof  Let $Z$ be a $(t+2)$-subspace of $V$. Observe that
\begin{align*}
\mathcal{H}_3(Z)=\left(\bigcup_{D\in{Z\brack t+1,0}}\left\{F\in{V\brack n,0}\mid F\cap Z=D\right\}\right)\cup\left\{F\in{V\brack n,0}\mid Z\subseteq F\right\}.
\end{align*}

Suppose the type of $Z$ is $(t+2,0)$. By Lemma~\ref{lem5}, the number of subspaces $F\in{V\brack n,0}$ with $Z\subseteq F$ is $N^\prime(t+2,0;n,0;n+\ell,n)=q^{\ell(n-t-2)}$, and for each $D\in{Z\brack t+1,0}$ the number of subspaces $F\in{V\brack n,0}$ satisfying $F\cap Z=D$ is $N^\prime(t+1,0;n,0;n+\ell,n)-N^\prime(t+2,0;n,0;n+\ell,n)=q^{\ell(n-t-1)}-q^{\ell(n-t-2)}.$ Since $\left|{Z\brack t+1,0}\right|=\theta_{t+2}$, we have (\ref{b-equ1}) holds.

Suppose the type of $Z$ is $(t+2,1)$. It is clear that the number of subspaces $F\in{V\brack n,0}$ with $Z\subseteq F$ is $0$. For each $D\in{Z\brack t+1,0}$, note that
\begin{align*}
\left\{F\in{V\brack n,0}\mid F\cap Z=D\right\}=\left\{F\in{V\brack n,0}\mid D\subseteq F\right\},
\end{align*}
implying that the number of subspaces $F\in{V\brack n,0}$ with $F\cap Z=D$ is $N^\prime(t+1,0;n,0;n+\ell,n)=q^{\ell(n-t-1)}.$ Since $\left|{Z\brack t+1,0}\right|=N^\prime(0,0;t+1,0;t+2,t+1)=q^{t+1}$, we have (\ref{b-equ10}) holds.

By (\ref{b-equ1}) and (\ref{b-equ10}), we have
\begin{align*}
	h_3(t+2,0)-h_3(t+2,1)=&q^{\ell(n-t-1)}\left(\theta_{t+2}-q^{t+1}\right)-q^{\ell(n-t-2)+1}\theta_{t+1}\\
	=&q^{\ell(n-t-1)}\theta_{t+1}-q^{\ell(n-t-2)+1}\theta_{t+1}>0,
\end{align*}
which implies that $h_3(t+2,0)>h_3(t+2,1).$
$\qed$

\begin{lemma}\label{b-lem14}
Suppose $4\leq n+1\leq \ell$ and $1\leq t\leq n-2.$ The following hold.
\begin{itemize}
\item[{\rm(i)}] If $1\leq t\leq\frac{n}{2}-\frac{3}{2}$, then $h_3(t+2,0)<f^\prime(n,\ell,t)$.
\item[{\rm(ii)}] If $t=\frac{n}{2}-1$, then $h_3(t+2,0)>f^\prime(n,\ell,t)>h_3(t+2,1)$.
\item[{\rm(iii)}] If $\frac{n}{2}-\frac{1}{2}\leq t\leq n-2$, then $h_3(t+2,1)>f^\prime(n,\ell,t)$.
\end{itemize}
\end{lemma}
\proof By (\ref{b-equ1}), we have
\begin{align*}
\frac{f^\prime(n,\ell,t)-h_3(t+2,0)}{q^{\ell(n-t-2)}\theta_{n-t-1}}
=q^{\ell}\left(q-\frac{q^{t+2}-1}{q^{n-t-1}-1}\right)-\frac{q^3(q^{n-t-2}-1)}{q^2-1}+\frac{q(q^{t+1}-1)}{q^{n-t-1}-1}.
\end{align*}
If $t\leq\frac{n}{2}-\frac{3}{2}$, then $q^{n-t-1}-1\geq q^{t+2}-1$, and
\begin{align*}
\frac{f^\prime(n,\ell,t)-h_3(t+2,0)}{q^{\ell(n-t-2)}\theta_{n-t-1}}>q^{\ell}-q^{n-t}>0
\end{align*}
by (\ref{b00}) and $\ell\geq n+t+1$. If $t=\frac{n}{2}-1$, then
\begin{align*}
\frac{f^\prime(n,\ell,t)-h_3(t+2,0)}{q^{\ell(n-t-2)}\theta_{n-t-1}}=\frac{-q^\ell(q-1)}{q^{t+1}-1}-\frac{q^3(q^t-1)}{q^2-1}+q<0.
\end{align*}

By (\ref{b-equ10}), we have
\begin{align}\label{b-ll-5}
\frac{f^\prime(n,\ell,t)-h_3(t+2,1)}{q^{\ell(n-t-2)}\theta_{n-t-1}}=& q^{\ell+1}\left(1-\frac{q^t(q-1)}{q^{n-t-1}-1}\right)-\frac{q^3(q^{n-t-2}-1)}{q^2-1}.
\end{align}
If $t=\frac{n}{2}-1$, by (\ref{b00}), then
\begin{align*}
\frac{f^\prime(n,\ell,t)-h_3(t+2,1)}{q^{\ell(n-t-2)}\theta_{n-t-1}}\geq \frac{q^{\ell+1}(q^t-1)}{q^{t+1}-1}-q^{t+2}>q^{\ell}-q^{t+2}>0.
\end{align*}
If $\frac{n}{2}-\frac{1}{2}\leq t\leq n-2$, then $q^{t}(q-1)>q^{n-t-1}-1$, and $f^\prime(n,\ell,t)<h_3(t+2,1)$ by (\ref{b-ll-5}). $\qed$
\begin{lemma}\label{b-lem13}
	Let $1\leq t\leq n-2$ and $n+t+2\leq \ell$ with $(\ell,q)\neq(n+t+2,2)$ or $(n+t+3,2).$ If $\mathcal{F}\subseteq{V\brack n,0}$ is a maximal $t$-intersecting family with $t+2\leq\tau_t(\mathcal{F})\leq n,$  then $|\mathcal{F}|<f^\prime(n,\ell,t)$.
\end{lemma}
\proof By Lemma~\ref{b-prop7} and (\ref{b00}), we have
\begin{align*}
	 \frac{f^\prime(n,\ell,t)}{|\mathcal{F}|}&\geq\frac{q^{\ell-2}(q^2-1)(q-1)^2}{(q^{n-t}-1)(q^{t+2}-1)(q^{t+1}-1)}-\frac{(q^{n-t-2}-1)(q-1)^2}{(q^{n-t}-1)(q^{t+2}-1)(q^{t+1}-1)}\\
	&>\frac{q^{\ell-2}(q^2-1)(q-1)^2}{q^{n-t}\cdot q^{t+2}\cdot q^{t+1}}-\frac{1}{q^{2}\cdot q^{t+1}\cdot q^{t}}=q^{\ell-n-t-5}(q^2-1)(q-1)^2-q^{-2t-3}.
\end{align*}

If $\ell\geq n+t+2$ and $q\geq 3$, then
\begin{align*}
	f^\prime(n,\ell,t)\cdot|\mathcal{F}|^{-1}>(1-q^{-2})(1-q^{-1})(q-1)-q^{-2t-3}>2(1-3^{-2})(1-3^{-1})-3^{-5}>1.
\end{align*}
If $\ell\geq n+t+4$ and $q=2$, then $f^\prime(n,\ell,t)/|\mathcal{F}|>1.5-2^{-5}>1,$ and the lemma holds.  $\qed$
\begin{lemma}\label{b-lem11}
Let $1\leq t\leq n-2$, and $n+t+2\leq \ell$ with $(\ell,q)\neq(n+t+2,2)$ or $(n+t+3,2).$ If $\mathcal{F}\subseteq{V\brack n,0}$ is a maximal $t$-intersecting family with $\tau_t(\mathcal{F})=t+1,$ and $\mathcal{F}$ is not a family given in Families~II and III, then $|\mathcal{F}|<f^\prime(n,\ell,t)$.
\end{lemma}
\proof Let $\mathcal{T}$ and $M$ be as in Hypothesis~\ref{b-hyp1}. Set
\begin{align*}
g(n+\ell,n,t)=q^{-\ell(n-t-2)}(\theta_{n-t-1})^{-1}(f^\prime(n,\ell,t)-|\mathcal{F}|).
\end{align*}

\noindent\textbf{Case 1.} $\mathcal{T}$ is not a $t$-intersecting family.

Observe that $\ell<3n-3t-2$ in this case from Lemma~\ref{b-lem2}.
\
\medskip
\

\noindent\textbf{Subcase 1.1.} $\ell\geq n+t+2$ and $q\geq 3$.

Then, we have $n-t\geq t+3\geq 4$, $n+t\geq 4+2t\geq 6$ and $\ell\geq 8.$ Observe that $$(q+1)^2(\theta_{n-t-1})^{-1}\leq (q+1)^2(\theta_{3})^{-1}=(q+1)^2(q^2+q+1)^{-1}<2,$$ $$(q^2+2q+2)\theta_{t-1}(\theta_{n-t-1})^{-1}\leq 3q^2\theta_{t-1}(\theta_{t+2})^{-1}<3q^{-1},$$ and $q^{-\ell+5}\theta_{t-2}<1$. By (\ref{b00}) Lemma~\ref{b-prop5}, we have
\begin{align*}
g(n+\ell,n,t)\geq&q^{\ell+1}-\frac{q^3(q^{n-t-2}-1)}{q^2-1}-\frac{(q+1)^2q^\ell}{\theta_{n-t-1}}-\frac{(q^2+2q+2)q^2\theta_{t-1}}{\theta_{n-t-1}}-\frac{q^{-\ell+5}\theta_{t-1}\theta_{t-2}}{\theta_{n-t-1}}   \\
\geq& q^{\ell+1}-q^{n-t}-2q^\ell-3q-1 \geq q^{\ell}-q^{n-t}-3q-1>0.	
\end{align*}

\noindent\textbf{Subcase 1.2.} $\ell\geq n+t+4$ and $q=2$.

Then $n-t\geq t+4\geq 5$, $n+t\geq 7$ and $\ell\geq 11$. Therefore
\begin{align*}
	&g(n+\ell,n,t)\\
	\geq&2^{\ell+1}-\frac{8(2^{n-t-2}-1)}{3}-\frac{9\cdot 2^\ell}{2^{n-t-1}-1}-\frac{40\cdot (2^{t-1}-1)}{2^{n-t-1}-1}-\frac{2^{-\ell+5}(2^{t-1}-1)(2^{t-2}-1)}{2^{n-t-1}-1}   \\
	\geq& 2^{\ell+1}-\frac{1}{3}\cdot 2^{n-t+1}+\frac{8}{3}-\frac{3}{5}\cdot 2^\ell-\frac{5}{2}-1
	>  \frac{7}{5}\cdot 2^{n+t+4}-\frac{1}{3}\cdot 2^{n-t+1}-\frac{5}{6}>0.	
\end{align*}	
\
\medskip
\

\noindent\textbf{Case 2.} $\mathcal{T}$ is a $t$-intersecting family.

By Lemmas~\ref{b-lem3}, \ref{b-prop3}, \ref{b-prop3-1} and \ref{b-prop4}, it suffices to show that $|\mathcal{F}|<f^\prime(n,\ell,t)$ when $\tau_t(\mathcal{T})=t$ and $t+1\leq\dim (M)\leq n-1.$
\
\medskip
\

\noindent\textbf{Case 2.1.} $\ell\geq n+t+2$ and $q\geq 3$.
\
\medskip
\

Observe that $\ell\geq n+t+2\geq\max\{n-t+4, 2t+4, n+3\}$, $q^3(\theta_{n-t-1})^{-1}{n-t-1\brack 2}<q^{n-t}$ and $q^{i-1}\leq\theta_{i}<q^{i}$ for $i\in\mathbb{Z}^{\rm +}$ by (\ref{b00}). From Lemma~\ref{b-prop1}, the following hold.

If $\dim(M)=t+1$, then
\begin{align*}
g(n+\ell,n,t)>&  q^{\ell+1}-q^{n-t}-q^\ell(\theta_{n-t-1})^{-1}-(q\theta_{n-t}-1)q^2\theta_{t+1}\\
\geq& q^{\ell+1}-q^{n-t}-q^{\ell-n+t+2}-q^3\theta_{n-t}\theta_{t+1}+q^2\theta_{t+1}  \\
\geq & 2q^{n+t+2}-q^{n-t}-q^{n+4}(q-1)^{-2}+q^2\theta_{t+1}> 0.
\end{align*}

If $\dim(M)=t+2\leq n-1$, then $\theta_{n-t-1}\geq q+1,$ and
\begin{align*}
	g(n+\ell,n,t)>& q^{\ell+1}-q^{n-t}-q^\ell-q^2(q^{n-t+1}+q^{t+2}-2q)  \\
	\geq & 2q^{\ell}-q^{n-t}-q^{n-t+3}-q^{t+4}+2q^3  \\
	\geq &  2q^{\ell}-q^{n-t+4}-q^{t+4}+2q^3> 0.
\end{align*}

If $t+3\leq \dim(M)\leq n-1$, then
\begin{align*}
	g(n+\ell,n,t)>&  q^{\ell+1}-q^{n-t}-q^\ell-q^3\theta_{n-t}-q^{t+3}(\theta_{n-t-1})^{-1}\\
\geq& 2q^\ell-q^{n-t+4}-q^{t+3}\geq 0.  \\
\end{align*}

\noindent\textbf{Case 2.2.} $\ell\geq n+t+4$ and $q=2$.
\
\medskip
\

Observe that $\ell\geq n+t+4\geq\max\{n-t+6, 2t+6, n+5\}$.
If $\dim(M)=t+1$, then
\begin{align*}
	g(n+\ell,n,t)\geq&  2^{\ell+1}-\frac{8}{3}\cdot(2^{n-t-2}-1)-2^\ell(2^{n-t-1}-1)^{-1}-8(2^{n-t}-1)(2^{t+1}-1)+4\theta_{t+1}\\
	\geq& 2^{\ell+1}-2^\ell-2^{n+4}+\frac{8}{3}+2^{n-t+3}-\frac{1}{3}\cdot 2^{n-t+1}+2^{t+4}-8+4\theta_{t+1}>0 .
\end{align*}

If $\dim(M)=t+2\leq n-1$, then $n-t\geq 3$ and
\begin{align*}
	g(n+\ell,n,t)>& 2^{\ell+1}-\frac{8}{3}\cdot(2^{n-t-2}-1)-3\cdot2^\ell(2^{n-t-1}-1)^{-1}-4(2^{n-t+1}+2^{t+2}-4) \\
	\geq& 2^{\ell+1}-2^\ell- \frac{13}{3}\cdot2^{n-t+1}-2^{t+4}+\frac{8}{3}+16 \\
	= &2^{\ell-1}-\frac{13}{3}\cdot2^{n-t+1}+2^{\ell-1}-2^{t+4}+\frac{56}{3}>0.
\end{align*}

If $t+3\leq \dim(M)\leq n-1$, then
\begin{align*}
	g(n+\ell,n,t)>& 2^{\ell+1}-\frac{8}{3}\cdot(2^{n-t-2}-1)-2^\ell-8(2^{n-t}-1)-2^{t+3}(2^{n-t-1}-1)^{-1}\\
	\geq& 2^{\ell-1}-\frac{13}{3}\cdot 2^{n-t+1}+2^\ell-2^{t+3}+\frac{32}{3}>0.
\end{align*}

Hence, the required result follows.   $\qed$

\
\medskip
\

\textbf{Proof of Theorem~\ref{b-main1}.}\quad

By Lemmas \ref{b-lem6} (i), \ref{b-lem14}, \ref{b-lem13} and \ref{b-lem11}, the desired result follows.  $\qed$
\section{The proof of Theorem~\ref{b-main2}}
In this section, we firstly give the exact value of the size of Family II, and then prove Theorem~\ref{b-main2} by comparing the size of each family given in Theorem~\ref{b-main1}.
\begin{lemma}\label{b-lem4}
Suppose $E$ is an $(n+r,r)$-subspace and $O$ is an $(m,0)$-subspace of $V$ with $O\subseteq E$. Then for each $a\in\{m,m+1,\ldots,n\}$ the number of $(a,0)$-subspaces $F$ of $V$ satisfying $F\cap E=O$ is
	\begin{align*}
	\prod_{j=0}^{a-m-1}\frac{(q^\ell-q^{r+j})(q^n-q^{m+j})}{q^a-q^{m+j}}.
\end{align*}
\end{lemma}
\proof Let $\mathcal{Z}(a)=\{F\in {V\brack a,0}\mid F\cap E=O\}$ for $a\in\{m,m+1,\ldots,n\}$. Observe that $|\mathcal{Z}(m)|=1.$ For each $a\in\{m,m+1,\ldots,n-1\}$, set
\begin{align*}
\mathcal{Z}=\{(F_1,F_2)\in\mathcal{Z}(a)\times \mathcal{Z}(a+1)\mid F_1\subseteq F_2\}.
\end{align*}

For a fixed $F_2\in\mathcal{Z}(a+1)$. Since $\{F_1\in \mathcal{Z}(a)\mid F_1\subseteq F_2\}=\{F_1\in{V\brack a,0}\mid O\subseteq F_1\subseteq F_2\}$, by Lemma~\ref{lem5}, the number of $F_1$ in $\mathcal{Z}(a)$ with $F_1\subseteq F_2$ is $\theta_{a+1-m}$, which implies that $|\mathcal{Z}|=|\mathcal{Z}(a+1)|\theta_{a+1-m}$. On the other hand, for a fixed $F_1\in\mathcal{Z}(a)$, Observe that
\begin{align*}
&\{F_2\in \mathcal{Z}(a+1)\mid F_1\subseteq F_2\}\\
=&\left\{F_2\in{V\brack a+1,0}\mid F_1\subseteq F_2\right\}\setminus	\left\{F_2\in{V\brack a+1,0}\mid F_1\subseteq F_2, \dim(F_2\cap E)=m+1\right\}\\
=&\left\{F_2\in{V\brack a+1,0}\mid F_1\subseteq F_2\right\}\setminus	\left\{F_2\in{V\brack a+1,0}\mid F_1\subseteq F_2\subseteq F_1+E\right\}.
\end{align*}
Since the type of $F_1+E$ is $(a+n+r-m, a+r-m)$, by Lemma~\ref{lem5}, we obtain that the number of $F_2$ in $\mathcal{Z}(a+1)$ with $F_1\subseteq F_2$ is $N^\prime(a,0;a+1,0;n+\ell,n)-N^\prime(a,0;a+1,0;a+n+r-m,n)=q^{a+r-m}(q^{\ell-a+m-r}-1)\theta_{n-a}$, which implies that $|\mathcal{Z}|=|\mathcal{Z}(a)|q^{a+r-m}(q^{\ell-a+m-r}-1)\theta_{n-a}.$ Hence, we have
\begin{align*}
|\mathcal{Z}(a)|q^{a+r-m}(q^{\ell-a+m-r}-1)\theta_{n-a}=|\mathcal{Z}(a+1)|\theta_{a+1-m}
\end{align*}
for $a\in\{m,m+1,\ldots,n-1\}$, and the result follows by induction.   $\qed$
\begin{lemma}\label{b-lem5}
	 Suppose $E$ is an $(n,1)$-subspace and $O$ is an $(m,0)$-subspace of $V$ with $O\subseteq E$. Then for each $a\in\{m+1,m+2,\ldots,n\}$ the number of $(a,0)$-subspaces $F$ of $V$ satisfying $F\cap E=O$ and $F\nsubseteq W+E$ is
	\begin{align*}
	q^{n+\ell-a}
	\cdot\prod_{j=1}^{a-m-1}\frac{(q^\ell-q^j)(q^n-q^{m+j})}{q^a-q^{m+j}}.
\end{align*}
\end{lemma}
\proof Let $\mathcal{Z}(a)=\{F\in {V\brack a,0}\mid F\cap E=O,\ F\nsubseteq W+E\}$ for $a\in\{m+1,m+2,\ldots,n\}$. Since $\mathcal{Z}(m+1)=\{F\in {V\brack m+1,0}\mid O\subseteq F\}\setminus\{F\in {V\brack m+1,0}\mid O\subseteq F\nsubseteq W+E\}$ and the type of $W+E$ is $(n+\ell-1,\ell)$, we have $|\mathcal{Z}(m+1)|=N^\prime(m,0;m+1,0;n+\ell,n)-N^\prime(m,0;m+1,0;n+\ell-1,n-1)=q^{n+\ell-m-1}$ by Lemma~\ref{lem5}. For each $a\in\{m+1,m+2,\ldots,n-1\}$, set
\begin{align*}
\mathcal{Z}=\{(F_1,F_2)\in\mathcal{Z}(a)\times \mathcal{Z}(a+1)\mid F_1\subseteq F_2\}.
\end{align*}

For a fixed $F_2\in\mathcal{Z}(a+1)$. Observe that
\begin{align*}
	&\{F_1\in \mathcal{Z}(a)\mid F_1\subseteq F_2\}\\
	=&\left\{F_1\in{V\brack a,0}\mid O\subseteq F_1\subseteq F_2,\ F_1\nsubseteq W+E\right\}\\
	=&\left\{F_1\in{V\brack a,0}\mid O\subseteq F_1\subseteq F_2\right\}\setminus \left\{F_1\in{V\brack a,0}\mid O\subseteq F_1\subseteq F_2\cap(W+E)\right\}.
\end{align*}
Since the type of $W+E$ is $(n+\ell-1,\ell)$ and $F_2\nsubseteq W+E,$ we obtain that the type of $F_2\cap(W+E)$ is $(a,0)$. Hence the number of $F_1\in \mathcal{Z}(a)$ with $F_1\subseteq F_2$ is $\theta_{a+1-m}-1=q\theta_{a-m}$, and  $|\mathcal{Z}|=|\mathcal{Z}(a+1)|q\theta_{a-m}$. On the other hand, for a fixed $F_1\in\mathcal{Z}(a)$, observe that
\begin{align*}
	&\{F_2\in \mathcal{Z}(a+1)\mid F_1\subseteq F_2\}\\
	=&\left\{F_2\in {V\brack a+1,0}\mid F_1\subseteq F_2, F_2\cap E=O\right\}\\
	=&\left\{F_2\in{V\brack a+1,0}\mid F_1\subseteq F_2\right\}\setminus	\left\{F_2\in{V\brack a+1,0}\mid F_1\subseteq F_2, \dim(F_2\cap E)=m+1\right\}\\
	=&\left\{F_2\in{V\brack a+1,0}\mid F_1\subseteq F_2\right\}\setminus	\left\{F_2\in{V\brack a+1,0}\mid F_1\subseteq F_2\subseteq F_1+E\right\}.
\end{align*}
Since the type of $F_1+E$ is $(a+n-m, a-m)$, we obtain that the number of $F_2$ in $\mathcal{Z}(a+1)$ with $F_1\subseteq F_2$ is $N^\prime(a,0;a+1,0;n+\ell,n)-N^\prime(a,0;a+1,0;a+n-m,n)=q^{a-m}(q^{\ell-a+m}-1)\theta_{n-a}$, which implies that $|\mathcal{Z}|=|\mathcal{Z}(a)|q^{a-m}(q^{\ell-a+m}-1)\theta_{n-a}.$ Hence, we have
\begin{align*}
|\mathcal{Z}(a)|q^{a-m}(q^{\ell-a+m}-1)\theta_{n-a}=|\mathcal{Z}(a+1)|q\theta_{a-m}
\end{align*}
for $a\in\{m+1,m+2,\ldots,n-1\}$, and the result follows by induction.   $\qed$

\begin{lemma}\label{b-lem1} The size of each family constructed in Family II only depends on the type of $X$, $M$ and $C$. Set $h_2(t,0;n,k;c,c-n)=|\mathcal{H}_2(X,M,C)|$ with $X\in{V\brack t,0}$, $M\in{V\brack n,k}$ and $C\in{V\brack c,c-n}$. Then the following hold:
	\begin{align}
		h_2(t,0;n,0;c,c-n)=&q^{\ell(n-t)}-\prod_{j=0}^{n-t-1}(q^\ell-q^j)+\prod_{i=0}^{c-n-1}(q^{n-t}-q^i)\nonumber\\
		&\cdot\prod_{j=0}^{2n-c-t-1}(q^\ell-q^{c-n+j})+q^{n-t}(q^{c-n}-1)\theta_{t},\label{b-equ3}\\
		h_2(t,0;n,1;c,c-n)=&q^{\ell(n-t)}-q^\ell\cdot\prod_{j=1}^{n-t-1}(q^\ell-q^j)+q^{n-t}\cdot\prod_{i=1}^{c-n-1}(q^{n-t}-q^i)\nonumber\\
		&\cdot\prod_{j=0}^{2n-c-t-1}(q^\ell-q^{c-n+j})+q^{c-t-1}(q^{t}-1).\label{b-equ4}
	\end{align}
In particular, set $h_1(t,0;n+1,1)=|\mathcal{H}_1(X,M)|$ with $X\in{V\brack t,0}$ and $M\in{V\brack n+1,1}$, and then
	\begin{align}
	h_1(t,0;n+1,1)=q^{\ell(n-t)}-\prod_{j=1}^{n-t}(q^\ell-q^j)+q^{n-t}(q^t-1).\label{b-equ2}
\end{align}
\end{lemma}
\proof Suppose $X\in{V\brack t,0}$, $M\in{V\brack n,k}$ and $C\in{V\brack c,c-n}$ are three subspaces with $X\subseteq M\subseteq C$, where $k\in\{0,1\}$ and $c\in\{n+1,n+2,\cdots,2n-t,n+\ell\}$. Let $\mathcal{A}(X,M,C)$, $\mathcal{B}(X,M,C)$ and $\mathcal{C}(X,M,C)$ be the families constructed in Family II. We have
\begin{align*}
	\mathcal{A}(X,M,C)=&\left\{F\in{V\brack n,0}\mid X\subseteq F\right\}\setminus\left\{F\in{V\brack n,0}\mid F\cap M=X\right\},\\
	\mathcal{B}(X,M,C)=&\bigcup_{D\in{C\brack c-n+t,0},\ D\cap M=X}\left\{F\in{V\brack n,0}\mid F\cap C=D\right\},\\
	\mathcal{C}(X,M,C)=&\bigcup_{D\in{M\brack n-1,0},\ X\nsubseteq D}\left\{F\in{V\brack n,0}\mid F\subseteq C,\ F\cap M=D\right\}.
\end{align*}

By Lemmas~\ref{lem5}, \ref{b-lem4} and \ref{b-lem5}, the size of $\left\{F\in{V\brack n,0}\mid X\subseteq F\right\}$ is $N^\prime(t,0;n,0;n+\ell,n)=q^{\ell(n-t)};$ and the number of subspaces $F$ in ${V\brack n,0}$ with $F\cap M=X$ is $\prod_{j=0}^{n-t-1}(q^\ell-q^j)$ if $k=0$, and is $q^\ell\cdot\prod_{j=1}^{n-t-1}(q^\ell-q^j)$ if $k=1$. Hence,
\begin{align*}
|\mathcal{A}(X,M,C)|=\left\{\begin{array}{ll}
q^{\ell(n-t)}-\prod_{j=0}^{n-t-1}(q^\ell-q^j),& \mbox{if}\ k=0,\\
q^{\ell(n-t)}-q^\ell\cdot\prod_{j=1}^{n-t-1}(q^\ell-q^j),& \mbox{if}\ k=1.\end{array}\right.
\end{align*}

For each $D\in{C\brack c-n+t,0}$ with $D\cap M=X$, if $D\in (C\cap W)+M$, then $\dim(D+M)=c-n+t+n-t=c\leq\dim((C\cap W)+M)$, implying that $C\cap W\cap M=0.$ Thus, we have $D\nsubseteq (C\cap W)+M$ for each $D\in{C\brack c-n+t,0}$ with $D\cap M=X$ if $M$ is of type $(n,1).$ By Lemmas \ref{b-lem4} and \ref{b-lem5}, the number of subspaces $D$ in ${C\brack c-n+t,0}$ with $D\cap M=X$ is $\prod_{i=0}^{c-n-1}(q^{n-t}-q^i)$ if $M$ is of type $(n,0)$, and is $q^{n-t}\cdot\prod_{i=1}^{c-n-1}(q^{n-t}-q^i)$ if $M$ is of type $(n,1)$. Observe that for each $D\in{C\brack c-n+t,0}$ the number of $(n,0)$-subspaces $F$ of $V$ satisfying $F\cap C=D$ is $\prod_{j=0}^{2n-c-t-1}(q^\ell-q^{c-n+j})$ by Lemma~\ref{b-lem4}. Hence,\begin{align*}
|\mathcal{B}(X,M,C)|=\left\{\begin{array}{ll}
\prod_{i=0}^{c-n-1}(q^{n-t}-q^i)\cdot\prod_{j=0}^{2n-c-t-1}(q^\ell-q^{c-n+j}),& \mbox{if}\ k=0,\\
q^{n-t}\cdot\prod_{i=1}^{c-n-1}(q^{n-t}-q^i)\cdot\prod_{j=0}^{2n-c-t-1}(q^\ell-q^{c-n+j}),& \mbox{if}\ k=1.\end{array}\right.
\end{align*}

By Lemma~\ref{lem5}, the number of $(n-1,0)$-subspaces $D$ of $M$ with $X\nsubseteq D$ is $N^\prime(0,0;n-1,0;n,n-k)-N^\prime(t,0;n-1,0;n,n-k)=q^{k(n-1)}{n-k\brack n-1}-q^{k(n-1-t)}{n-k-t\brack n-1-t}$, and for each $D\in{M\brack n-1,0}$ with $X\nsubseteq D$ the number of $(n,0)$-subspaces $F$ of $V$ with $F\subseteq C$ and $F\cap M=D$ is $N^\prime(n-1,0;n,0;c,n)-N^\prime(n,k;n,0;c,n)=q^{c-n}-q^{k(c-n)}{c-n-k\brack -k}$. Hence,
\begin{align*}
|\mathcal{C}(X,M,C)|=\left\{\begin{array}{ll}
q^{n-t}(q^{c-n}-1)\theta_{t},& \mbox{if}\ k=0,\\
q^{c-t-1}(q^{t}-1),& \mbox{if}\ k=1.
\end{array}\right.
\end{align*}

Therefore, (\ref{b-equ3}) and (\ref{b-equ4}) hold, and (\ref{b-equ2}) holds by Remark~\ref{b-obs}.  $\qed$
\begin{lemma}\label{b-lem15}
	Suppose $4\leq n+1\leq \ell$, $1\leq t\leq n-2$ and $k\in\{0,1\}$. Then the following hold.
	\begin{itemize}
		\item[{\rm(i)}] We have
		\begin{align*}
			h_1(t,0;n+1,1)=&h_2(t,0;n,k;n+1,1)>h_2(t,0;n,k;n+2,2)>\cdots\\
			&>h_2(t,0;n,k;2n-t,n-t).
		\end{align*}
		\item[{\rm(ii)}] We have $h_2(t,0;n,0;n+\ell,\ell)>h_2(t,0;n,1;n+\ell,\ell)$.
		\item[{\rm(iii)}]If $1\leq t\leq n-3$, then $h_1(t,0;n+1,1)>h_2(t,0;n,0;n+\ell,\ell).$
		\item[{\rm(iv)}] If $t=n-2$, then $h_2(t,0;n,0;n+\ell,\ell)>h_1(t,0;n+1,1)$.
	\end{itemize}
\end{lemma}
\proof (i)\quad Notice that $h_1(t,0;n+1,1)=h_2(t,0;n,0;n+1,1)=h_2(t,0;n,1;n+1,1)$ by Remark~\ref{b-obs}. For each $c\in\{n+1,n+2,\ldots,2n-t-1\}$, from $\ell\geq n+1$, we have
\begin{align*}
	&h_2(t,0;n,0;c,c-n)-h_2(t,0;n,0;c+1,c+1-n)\\
	=&(q^\ell-q^{n-t})\cdot \prod_{i=0}^{c-n-1}(q^{n-t}-q^i)\cdot\prod_{j=1}^{2n-c-t-1}(q^{\ell}-q^{c-n+j})-q^{c-t}(q^t-1)\\
	\geq&(q^\ell-q^{n-t})(q^{n-t}-1)-q^{c-t}(q^t-1)\geq q^{n-t}(q^{t+1}-1)(q^{n-t}-1)-q^{2n-2t-1}(q^t-1)>0,
\end{align*}
and
\begin{align*}
	&h_2(t,0;n,1;c,c-n)-h_2(t,0;n,1;c+1,c+1-n)\\
	=&q^{n-t}(q^\ell-q^{n-t})\cdot \prod_{i=1}^{c-n-1}(q^{n-t}-q^i)\cdot\prod_{j=1}^{2n-c-t-1}(q^{\ell}-q^{c-n+j})-q^{c-t-1}(q-1)(q^t-1)\\
	>&q^{n-t+\ell-1}-q^{c}\geq q^{n-t+\ell-1}-q^{2n-t-1}>0.
\end{align*}
 Hence, the function $h_2(t,0;n,k;c,c-n)$ is decreasing as $c\in\{n+1,n+2,\ldots,2n-t\}$ increases, and (i) holds.
\
\medskip
\

(ii)\quad By (\ref{b-equ3}) and (\ref{b-equ4}), we have
\begin{align*}
	h_2(t,0;n,0;n+\ell,\ell)-h_2(t,0;n,1;n+\ell,\ell)=\prod_{j=1}^{n-t-1}(q^\ell-q^j)+q^{n-t}(q^t-1)\theta_{\ell-1}>0,
\end{align*}
implying that (ii) holds.
\
\medskip
\

(iii)\quad From Lemma~\ref{b-lem1}, $t\leq n-3$ and $n+1\leq \ell$, we have
\begin{align*}
	h_1(t,0;n+1,1)-h_2(t,0;n,0;n+\ell;\ell)=&(q^{n-t}-1)\cdot\prod_{j=1}^{n-t-1}(q^\ell-q^j)-q^{n-t+1}(q^t-1)\theta_{\ell-1}\\
	>&(q^{n-t}-1)(q^\ell-q)(q^\ell-q^2)-q^{n+\ell}\\
	>&q^{2\ell+n-t-3}-q^{n+\ell}>0.
\end{align*}
 Hence (iii) holds.
\
\medskip
\

(iv)\quad If $t=n-2$, by Lemma~\ref{b-lem1}, then
\begin{align*}
	h_2(t,0;n,0;n+\ell,\ell)-h_1(t,0;n+1,1)=&-(q^{\ell}-q)(q^2-1)+q^3(q^{n-2}-1)\theta_{\ell-1}\\
	>&(q^{\ell}-q)\left(q^2\theta_{n-2}-q^2+1\right)>0,
\end{align*}
which implies that (iv) holds.  $\qed$

\begin{lemma}\label{b-lem16}
	Suppose $4\leq n+1\leq \ell$ and $1\leq t\leq n-2.$ The following hold.
	\begin{itemize}
		\item[{\rm(i)}] If $t=\frac{n}{2}-1$, then $h_1(t,0; n+1,1)>h_3(t+2,0)$.
		\item[{\rm(ii)}] If $\frac{n}{2}-\frac{1}{2}\leq t\leq n-3$, then $h_3(t+2,0)>h_1(t,0; n+1,1)$.
		\item[{\rm(iii)}]If $t=n-2$, then $h_3(t+2,0)=h_2(t,0;n,0;n+\ell,\ell)$.
	\end{itemize}
\end{lemma}
\proof (i) Let $X$ be a $(t,0)$-subspace and $M$ be an $(n+1,1)$-subspace of $V$ with $X\subseteq M$. Suppose
\begin{align*}
	\mathcal{A}_{t+1}(X,M)=\bigcup_{X\subseteq D,\ D\in{M\brack t+1,0}}\left\{F\in{V\brack n,0}\mid F\cap M=D\right\}.
\end{align*}
Since the number of $(t+1,0)$-subspace $D$ of $M$ with $X\subseteq D$ is $N^\prime(t,0;t+1,0;n+1,n)=q\theta_{n-t}$, by Lemma~\ref{b-lem4}, we have
\begin{align*}
|\mathcal{A}_{t+1}(X,M)|=q\theta_{n-t}\cdot\prod_{j=1}^{n-t-1}(q^\ell-q^j).
\end{align*}

Since $t=\frac{n}{2}-1$, by (\ref{b-equ1}) and the construction of Family I, we have
\begin{align*}
h_1(t,0; n+1,1)-h_3(t+2,0)>q\theta_{t+2}\cdot\prod_{j=1}^{t+1}(q^\ell-q^j)-q^{\ell(t+1)}\theta_{t+2}+q^{\ell t+1}\theta_{t+1}.
\end{align*}
Suppose $q\cdot\prod_{j=1}^{t+1}(q^\ell-q^j)=q^{\ell(t+1)+1}+g(q)$, where the polynomial $g(q)$ on $q$ have at most $2^{t+1}-1$ terms and $\deg(g(q))=\ell t+t+2.$ Since $\ell\geq n+1\geq 2t+3$, we have
\begin{align*}
q\cdot\prod_{j=1}^{t+1}(q^\ell-q^j)-q^{\ell(t+1)}=(q-1)\cdot q^{\ell(t+1)}-g(q)=(q-1)\cdot q^{\ell-t-2}\cdot q^{\ell t+t+2}-g(q)>0,
\end{align*}
which implies that $h_1(t,0; n+1,1)>h_3(t+2,0)$.
\
\medskip
\

(ii)\quad Suppose $\frac{n}{2}-\frac{1}{2}\leq t\leq n-3$. From $n\leq 2t+1$ and $\ell\geq n+1\geq t+4$, observe that
\begin{align*}
	\frac{h_3(t+2,0)-q^{\ell(n-t-1)+1}\theta_{n-t}}{q^{\ell(n-t-2)}}=& q^\ell\theta_{t+2}-q\theta_{t+1}-q^{\ell+1}\theta_{n-t}\\
	=&(q-1)^{-1}\cdot\left(q^\ell(q^{t+2}-q^{n-t+1}+q-1)-q(q^{t+1}-1)\right)>0.
\end{align*}
 By Lemma~\ref{b-lem6}, we have $h_3(t+2,0)>h_1(t,0; n+1,1)$ as desired.
\
\medskip
\

(iii)\quad By Remark~\ref{b-obs}, the desired result follows.  $\qed$

\
\medskip
\

\textbf{Proof of Theorem~\ref{b-main2}.}\quad

By Theorem~\ref{b-main1}, Remark~\ref{b-obs}, Lemmas~\ref{b-lem15} and \ref{b-lem16}, the desired result follows. $\qed$

\section*{Acknowledgement}

B. Lv is supported by NSFC (12071039), K. Wang is supported by the National Key R\&D Program of China (No. 2020YFA0712900) and NSFC (12071039).

%\end{CJK*}

\end{document}